\documentclass[11pt,a4wide,oneside,reqno]{amsart}
\usepackage{amssymb,amsmath,url,a4wide}
\usepackage[utf8]{inputenc}
\usepackage{header}
\usepackage[left=3cm,right=3cm,top=3cm,bottom=2cm]{geometry}
\usepackage{enumitem}

\newcommand{\pr}{\mathrm{pr}}
\newcommand{\conf}{\mathcal{C}}
\newcommand{\niconf}{\mathcal{C^\#}}
\newcommand{\iconf}{\mathcal{C}^\mathrm{inj}}
\newcommand{\g}{\mathfrak{D}}
\newcommand{\citep}[1]{\cite{#1}}

\newcommand{\AQ}{{\Phi_{\text{AQ}}}}
\newcommand{\AH}{{\Phi_{\text{AH}}}}

\setlist[enumerate]{label=(\roman*)}
\numberwithin{equation}{section}
\usepackage[url=false,doi=true,maxnames=4]{biblatex}
\usepackage{xpatch}
\xpatchbibmacro{volume+number+eid}{%
  \setunit*{\adddot}%
}{%
}{}{}
\DeclareFieldFormat[article]{number}{\mkbibparens{#1}}
\renewbibmacro{in:}{%
  \ifentrytype{article}{}{\printtext{\bibstring{in}\intitlepunct}}}
\title[Linear configurations containing 4-APs are uncommon]{Linear configurations containing 4-term arithmetic progressions are uncommon}
\author{Leo Versteegen}
\address{Department of Pure Mathematics and Mathematical Statistics, Centre for
Mathematical Sciences, Wilberforce Road, Cambridge CB3 0WB, United Kingdom}
\email{lvv23@dpmms.cam.ac.uk}
\begin{document}
\begin{abstract}
\noindent A linear configuration is said to be common in $G$ if every 2-coloring of $G$ yields at least the number of monochromatic instances of a randomly chosen coloring. Saad and Wolf asked whether, analogously to a result by Thomason in graph theory, every configuration containing a 4-term arithmetic progression is uncommon. We prove this in $\F_p^n$ for $p\geq 5$ and large $n$ and in $\Z_p$ for large primes $p$.
\end{abstract}
\maketitle
\section{Introduction}
In \cite{goodman}, Goodman determined the minimal number of monochromatic triangles in any 2-coloring of the edges of a complete graph. Later, Erdős  \cite{erdos-common-conjecture} pointed out that this minimum, roughly $\binom{n}{3}/4$ for a complete graph on $n$ vertices, is (at least asymptotically) achieved by a random coloring. This motivated the following definition. Let $H$ be a graph on $v(H)$ vertices with $e(H)$ edges. For a 2-coloring $\chi\colon [n]^{(2)}\rightarrow \{1,2\}$ of the edges the complete graph $K_n$, let $m_H(\chi)$ be the number of monochromatic copies of $H$ with respect to the coloring $\chi$. The graph $H$ is called \textit{common} if
\begin{align}\label{eq:graph-common}
\liminf_{n\rightarrow \infty} \frac{\min \{ m_H(\chi) : \chi \colon [n]^{(2)}\rightarrow \{1,2\} \}}{\binom{n}{v(H)}} \geq 2^{-e(H)+1}.
\end{align}
Erdős conjectured that all complete graphs were common, but this was disproved by Thomason, who showed in \cite{thomason-erdos-disproof} that even the complete graph on four vertices is uncommon. Later on, Jagger, Šťovíček and Thomason \cite{containing-k4-uncommon} proved the following vastly more general result.
\begin{theorem}[Jagger, Šťovíček and Thomason]\label{thm:containing-k4-uncommon}
Every graph containing $K_4$ as a subgraph is uncommon.
\end{theorem}
This theorem refutes Erdős's conjecture quite robustly. On the other hand, several non-trivial classes of common graphs have been established (see \cite{thomason-graph-products-mono-mult,conlon-ramsey-mult-complete-graphs}, for example), making the study of common graphs a challenging problem of ongoing interest.\\\\
In \cite{saad-wolf-sidorenko-common}, Saad and Wolf introduced an arithmetic version of the notion of commonness, motivated in part by the wish to gain insights that might be transferred back to the graph-theoretic problem. In the arithmetic setting, the elements of a finite Abelian group are 2-colored, and the number of monochromatic instances of a given \textit{linear configuration} is compared to the expected number of monochromatic instances arising from a random coloring.\\\\
Even before \cite{saad-wolf-sidorenko-common}, Wolf had proved in \cite{wolf-4ap-uncommon} that the 4-term arithmetic progression (4-AP) is uncommon, using a then unpublished construction by Gowers, which is now available as \cite{gowers-4ap-uniform-non-sidorenko}. This is in contrast to the 3-AP, which is rather easily seen to be common, resembling the graph theoretic setting. In this light, the following question of Saad and Wolf about an arithmetic analogue of \Cref{thm:containing-k4-uncommon} is natural.
\begin{question}[Saad and Wolf]\label{qu:ap4-uncommon}
Is it true that every linear configuration containing a 4-AP is uncommon?
\end{question}
The main purpose of this paper is to give an affirmative answer to this question in the groups $\F_p^n$ for $p\geq 5$ and large $n$ and in $\Z_p$ for large primes $p$.\\\\
Throughout, we work in a finite Abelian group $G$, and we consider for $r\leq d\in \N$ a \textit{system of forms} $\Phi=(\phi_1,\ldots,\phi_d)$ where each $\phi_i$ is a map of the form
\begin{align*}
\phi_i\colon G^r \rightarrow G \qquad w\mapsto v_i^Tw
\end{align*}
for some $v_i\in \Z^r$. We may also view $\Phi$ as a map $G^r\rightarrow G^d$ and consider its image $\conf(\Phi):=\im \Phi$. We call any set $\conf\subseteq G^d$ that is given by the image of a system of forms a \textit{(linear) configuration}. The elements of $\conf$ are called the \textit{instances} of the configuration. Although the answers to many of the questions we may ask about linear configurations depend on the group they lie in, a system of forms may be encoded independently of the group $G$ by the vectors $v_i\in \Z^r$ or, equivalently, by a matrix $M\in \Z^{d\times r}$. We will say $M$ \textit{induces} $\Phi$ in $G$ if for all $i\in [d]$ and $w\in G^r$ it holds that $(Mw)_i=M_{i1}w_1+\ldots+M_{ir}w_r=\phi_i(w)$.\\\\
Linear configurations can equivalently be defined as the solution set to a system of homogeneous linear equations $Lx=0$ for $L\in \Z^{(d-r)\times d}$. In fact, this is the principal viewpoint in \cite{saad-wolf-sidorenko-common}. Both have their merit, and by looking at $\conf$ simply as a subset of $G^d$ in the following definition of commonness, we remain deliberately agnostic as to which viewpoint is taken.
\begin{definition}(Arithmetic multiplicity and common configurations)
For a linear configuration $\conf \subseteq G^d$ and $f\colon G\rightarrow \R$, define the \textit{arithmetic multiplicity} as\footnote{Here and throughout the paper we write for a set $X$ and a function $f\colon X\rightarrow \C$ the \textit{average of $f$ over $X$} as $\E(f)=\E_{x\in X} f(x)=\frac{1}{\lv X\rv}\sum_{x\in X} f(x)$.}
\begin{align*}
t_\conf(f):=\E_{v\in \conf} \prod_{i=1}^d f(v_i)=\frac{1}{\vert \conf\vert} \sum_{v\in \conf} \prod_{i=1}^d f(v_i).
\end{align*}
For $A\subseteq G$, we write $t_\conf(A)$ to mean $t_\conf(1_A)$, where $1_A$ is the indicator function of $A$. A linear configuration $\conf\subseteq G^d$ is called \textit{common} if for all $A\subseteq G$
\begin{align*}
t_{\conf}(A) + t_{\conf}(A^C) \geq \lp\frac{1}{2}\rp^{k-1}.
\end{align*}
We call the configuration \textit{fully common} if for all $f\colon G\rightarrow [0,1]$ 
\begin{align*}
t_\conf(f)+t_\conf(1-f)\geq \lp\frac{1}{2}\rp^{k-1}.
\end{align*}
\end{definition}
Our definition of commonness differs slightly from that in \cite{saad-wolf-sidorenko-common}, which is defined asymptotically so as to resemble the graph-theoretic notion given by \eqref{eq:graph-common} more closely. In the context of this paper, choosing between the definitions is purely a matter of taste and our results may easily be reformulated using the definition of commonness from \cite{saad-wolf-sidorenko-common} with only minor changes to the proofs.\\\\
It remains to set out what it means for a configuration to contain a 4-AP. It is possible to define in great generality when a configuration contains another one, but for the purpose of answering \Cref{qu:ap4-uncommon} we instead say very concretely that a configuration given by a system of forms $\Phi=(\phi_i)_{i\in [d]}$ in $r$ variables \textit{contains} a 4-AP if there is a set $\{i_1,i_2,i_3,i_4\}\in [d]^{(4)}$ such that 
\begin{align*}
\forall w\in G^r\colon \phi_{i_4}(w)-\phi_{i_3}(w)=\phi_{i_3}(w)-\phi_{i_2}(w)=\phi_{i_2}(w)-\phi_{i_1}(w).
\end{align*}
We are now ready to state our main results.
\begin{restatable}{theorem}{uncommonapvector}\label{thm:ap4vector}
Let $p\geq 5$ be a prime and let $M\in \Z^{d\times r}$. There exists $\e>0$ such that for all sufficiently large $n\in \N$ the following holds. If the system of forms $\Phi$ induced by $M$ on $G=\F_p^{n+1}$ is injective, contains a 4-AP and consists of pairwise distinct forms, then there exists a function $f\colon G\rightarrow [0,1]$ with
\begin{align}\label{eq:goal-Fpn}
t_{\Phi}(f)+t_{\Phi}(1-f) < \lp \frac{1}{2}\rp^{d-1} - \e.
\end{align}
\end{restatable}
The requirement that the forms be distinct is necessary, as is easily seen by taking one of the available lower bounds on the number of monochromatic 4-APs in a 2-coloring of $\F_p^n$ or $\Z_p$ (see \cite{lu-peng-4ap}, for example) and adding copies of one of the forms until the configuration becomes common. Note that the forms constituting the 4-AP cannot be distinct if $p\in \{2,3\}$ so the condition that $p\geq 5$ is void. Because the error parameter $\e$ in \Cref{thm:ap4vector} does not depend on $n$, we are able to obtain the following corollary.
\begin{restatable}{corollary}{uncommonapvectorcorollary}\label{cor:ap4vector}
Let $p\geq 5$ be a prime and let $M\in \Z^{d\times r}$. For all sufficiently large $n\in \N$ the following holds. If the system of forms $\Phi$ induced by $M$ on $G=\F_p^{n+1}$ is injective, contains a 4-AP and consists of distinct forms, then $\Phi$ is uncommon in $G$.
\end{restatable}
For the cyclic group $\Z_p$, we obtain the following result.
\begin{restatable}{theorem}{uncommonapcyclic}\label{thm:ap4cyclic}
Let $M\in \Z^{d\times r}$. There exists $\e>0$ such that for all sufficiently large primes $p$ the following holds. If the system of forms $\Phi$ induced by $M$ on $G=\Z_p$ is injective, contains a 4-AP and consists of distinct forms, then there exists a function $f\colon G\rightarrow [0,1]$ with
\begin{align}\label{eq:goal-Zp}
t_{\Phi}(f)+t_{\Phi}(1-f) < \lp \frac{1}{2}\rp^{d-1} - \e.
\end{align}
\end{restatable}
This theorem also has a matching corollary.
\begin{restatable}{corollary}{uncommonapcycliccorollary}\label{cor:ap4cyclic}
Let $M\in \Z^{d\times r}$. For all sufficiently large primes $p$ the following holds. If the system of forms $\Phi$ induced by $M$ on $G=\Z_p$ is injective, contains a 4-AP and consists of distinct forms, then $\Phi$ is uncommon in $G$.
\end{restatable}
Our proofs build on the constructions of Gowers \cite{gowers-4ap-uniform-non-sidorenko} and Wolf \cite{wolf-4ap-uncommon} mentioned earlier. In fact, they may be seen as an abstraction of their ideas resulting in a general strategy for building functions that yield low arithmetic multiplicities, which we call the ``muting method''. \Cref{thm:ap4vector} and \Cref{thm:ap4cyclic} should be thought of as exemplary applications of the muting method. Indeed, with only very minor changes to the proofs one can show, for example, that all configurations containing a 5-AP with the middle point missing, i.e., $(x-2y,x-y,x+y,x+2y)$, are uncommon.\\\\
A full characterization of common configurations remains elusive. To date, such a characterization exists only for configurations given by a single equation (or equivalently by systems of forms in $d-1$ variables). Specifically, Cameron, Cilleruelo and Serra \cite{odd-variable-equations-common} proved that if the equation has an odd number of variables (meaning that $d$ is odd), then the configuration is always common, even in non-Abelian groups. Saad and Wolf \cite{saad-wolf-sidorenko-common} proved that if $G=\F_p^n$ and $d$ is even, an equation is common if the coefficients of the equation can be partitioned into pairs such the coefficients in each pair sum to zero.\footnote{They in fact proved the stronger statement that configurations which satisfy this condition are \textit{Sidorenko}, see \cite{saad-wolf-sidorenko-common}.} They also conjectured that the existence of such a partition is a necessary condition for the equation to be common. Fox, Pham and Zhao confirmed this conjecture in \cite{fox-pham-zhao-sidorenko-equations}, and the author generalized their result from $\F_p^n$ to all finite Abelian groups in which the configuration is non-degenerate \cite{versteegen-sidorenko-common-abelian}.\\\\
Shortly after the present preprint was posted on the arXiv, the author learned that \Cref{thm:ap4vector} has been proved independently by Kamčev, Liebenau and Morrison \cite{kamcev-liebenau-morrison-uncommon-systems}. Additionally, they proved the more general result that in $\F_p^n$, any configuration that can be reparametrized such that it contains four forms depending only on two variables is uncommon. However, it seems that their method of proof does not extend to $\Z_p$.\\\\
Much remains to be understood about common configurations. We conclude the article by presenting the following example.
\begin{restatable}{theorem}{cubemissingcorner}\label{thm:cube-missing-corner}
Let $G$ be any Abelian group. The 3-cube with one missing vertex, given by the system $\Phi\colon G^4 \rightarrow G^7$ with
\begin{align*}
\Phi(x,h_1,h_2,h_3)=(x,x+h_1,x+h_2,x+h_3,x+h_1+h_2,x+h_1+h_3,x+h_2+h_3)
\end{align*}
is fully common in $G$.
\end{restatable}
This example shows that a configuration may be common for ``non-trivial'' reasons, meaning that it is neither given by a single equation nor does it have a global symmetry that allows for a convexity argument, as is the case for the full 3-cube. To see that this configuration is common it will be necessary to carefully dissect how different subconfigurations affect the overall arithmetic multiplicity. This illustrates the limitations of the muting method and highlights one of the obstacles to finding a characterization of common configurations.\\\\
We shall now lay out out how the rest of the paper is structured. In \Cref{sec:preliminaries}, we introduce a number of notions and tools that are going to be used throughout the article. \Cref{sec:4ap-vector} introduces the muting method in $\F_p^n$ and contains the proof of \Cref{thm:ap4vector}. In \Cref{sec:4ap-cyclic}, we adapt the muting method to the cyclic group $\Z_p$ and prove \Cref{thm:ap4cyclic}. We strongly encourage the reader to read \Cref{sec:4ap-vector} before \Cref{sec:4ap-cyclic} as the treatment of cyclic groups relies on the tools developed in the preceding section. Lastly, in \Cref{sec:cube-missing-corner} we prove \Cref{thm:cube-missing-corner} and discuss the limitations of the muting method.
\subsection*{Acknowledgments}
The author is grateful to be funded by Trinity College of the University of Cambridge through the Trinity External Researcher Studentship. Further, he wishes to express his deep gratitude to his master thesis supervisor Mathias Schacht  for introducing him to the area of additive combinatorics, to Emil Powierski for many hours of enlightening discussion on the topic and to his PhD supervisor  Julia Wolf for her generous help in improving this paper.
\section{Preliminiaries}\label{sec:preliminaries}
Our notation follows for the most part that of \cite{saad-wolf-sidorenko-common} and all symbols are standard unless explicitly introduced. For the avoidance of doubt, given a natural number $k$, we write $[k]$ for the set $\{1,2,\ldots, k\}$ and $[k]_0$ for the set $[k]\cup \{0\}$. Given any set $X$, we denote by $X^{(k)}$ the set of $k$-element subsets of $X$.\\\\
\subsection{Reparametrization and $C$-fractions}
Let $\Phi=(\phi_i)_{i\in [d]}$ be a system of forms describing a linear configuration. The proofs of both \Cref{thm:ap4vector} and \Cref{thm:ap4cyclic} rely on a careful analysis of small \textit{subconfigurations}, by which we mean subfamilies $(\phi_{i})_{i\in I}$ for some $I\subseteq [d]$. In what follows, we will prove a lemma that allows us to reparametrize systems of forms. This enables us to analyze the contributions of subconfigurations locally, i.e., without having to worry about how they fit into the global structure of the original configuration. When dealing with the cyclic group $\Z_p$, it will be important for us to keep tight control over the coefficients of the reparametrized forms. This control is expressed through the notion of $C$-fractions.
\begin{definition}[$C$-fractions]\label{def:C-fraction}
Let $C\in \N$ and $p$ a prime. An element $x\in \Z$ is called a \textit{$C$-fraction mod $p$} if there exist $y\in \Z, z\in \N$ such that $\vert y\vert,\vert z\vert < C$ and $xz\equiv y$ mod $p$. Typically, $p$ will be clear from the context and so we will simply call $x$ a $C$-fraction. We will also call $x\in \Z_p$ a $C$-fraction if it is the residue class of a $C$-fraction in $\Z$. We say that $C$ is a \textit{fraction bound} for $x$.
\end{definition}
We record some basic properties of $C$-fractions.
\begin{lemma}[Basic properties of $C$-fractions]\label{lem:C-fraction-properties}~
\begin{enumerate}
\item The number of $C$-fractions in $\Z_p$ is bounded by $4C^2$.
\item If $x_1$ is a $C$-fraction, then so is every $x_2$ with $x_1\equiv x_2 \mod p$. 
\item Let $x_1$ be a $C_1$-fraction and let $x_2$ be a $C_2$-fraction. Then $x_1x_2$ a is a $C_1C_2$-fraction, $x_1+x_2$ is a $2C_1C_2$-fraction, $-x_1$ is a $C_1$-fraction and any representative $x_1\inv$ of the inverse of $x_1$ mod $p$ is a $C_1$-fraction.
\end{enumerate}
\end{lemma}
\begin{proof}
The first two statements follow directly from the definition. The same holds for the part of (iii) concerning additive and multiplicative inverses. For the rest, let $y_1,y_2,z_1,z_2$ be such that $x_iz_i\equiv y_i$ for $i \in [2]$. Then
\begin{align*}
(x_1+x_2)z_1z_2&\equiv y_1z_2+y_2z_1 \qquad \mod p &&\text{and}& x_1x_2z_1z_2&\equiv y_1y_2\qquad \mod p,
\end{align*}
completing the proof.
\end{proof}
We now prove our reparametrization lemma, which is applicable both in $\F_p	n$ and $\Z_p$.
\begin{lemma}[Reparametrization Lemma]\label{lemma:reparametrization}
Let $\Phi$ be a system of $d$ forms in $r$ variables over a group $G=\F_p^n$ (where $n=1$ is possible), let $k\in [r]$ and let $\{i_1,\ldots,i_k\}\subseteq [d]$ be a set of indices such that $\phi_{i_1},\ldots,\phi_{i_k}$ are linearly independent $\mod p$. Then there exists a system $\Phi'$ with $\phi_{i_j}'(v)=v_j$ for $j\in [k]$ and $\conf(\Phi')=\conf(\Phi)$. Furthermore, if $\Phi$ is induced by a matrix whose coefficients are bounded by $C_1$, then $\Phi'$ is induced by a matrix all of whose coefficients are $C_2$-fractions for $C_2=rC_1(C_1\sqrt{r})^r$.
\end{lemma}
\begin{proof}
For $i\in [d]$, let $\xi_i\in \Z^r$ be such that the coefficients of $\xi_i$ are bounded by $C_1$ and $\phi_i(v)=\xi_i^Tv$. The family $(\xi_{i_j})_{j\in [k]}$ can be extended by adding $r-k$ standard unit vectors $\psi_1,\ldots,\psi_{r-k}$ without adding a linear dependence mod $p$. Consider the matrix
\begin{align*}
B=\begin{pmatrix}
\xi_{i_1}&\cdots&\xi_{i_k}&\psi_1&\cdots&\psi_{r-k}
\end{pmatrix}^T\in \Z^{r\times r}.
\end{align*}
There exists a matrix $A\in \Z^{r\times r}$ (the adjugate of $B$) such that $AB=\det(B)\cdot I_{r\times r}$; because the rows of $B$ are independent mod $p$ there exists $m\in \Z$ such that $m\det(B)\equiv 1$ mod $p$. Let now $M=mA$, $\phi_i'(v)=\xi_i^TMv$, and $\Phi'=(\phi_1',\ldots,\phi_d')$. For all $j\in [k]$ and $v\in G^r$ we have
\begin{align*}
\phi_{i_j}'(v)=m\cdot\xi_{i_j}^T Av=m\cdot(BA)(v)_j=m\det(B)v_j=v_j.
\end{align*}
Furthermore, note that the matrix $M$ acts bijectively on $G^r$ as it has an inverse given by $B$. Therefore,
\begin{align*}
\conf(\Phi')=\{\Phi(Mv)\colon v\in G^r\}=\{\Phi(v)\colon v\in G^r\}=\conf(\Phi).
\end{align*}
The system $\Phi'$ is induced by the matrix $N\in \Z^{d\times r}$ with $N_{ij}=(\xi_i^TM)_j=m(\xi_i^TA)_j$. By definition of $m$, we have $\det(B)N_{ij}\equiv (\xi_i^TA)_j$ mod $p$ so to complete the proof it suffices to show that $\det(B)$ and $(\xi_i^TA)_j$ are bounded by $C_2$. Recall that $C_1$ is an upper bound for all the coefficients occurring in the $\xi_i$ and hence also for those in $B$. By the definition of the adjugate matrix and, for example, Hadamard's inequality, we have that both the coefficients of $A$ and $\det(B)$ are bounded by $(C_1\sqrt{r})^r$. Thus $\lv \det(B)\rv <C_2$ and 
\begin{align*}
\lv (\xi_i^TA)_j\rv \leq m \sum_{k=1}^r \lv (\xi_i)_k A_{kj}\rv \leq rC_1(C_1\sqrt{r})^r=C_2,
\end{align*}
as desired.
\end{proof}
\subsection{Converting maps to sets}
When showing that a configuration is not common, it is often more convenient to construct a counterexample as a function $f\colon G\rightarrow[0,1]$ and convert it to a subset of $G$ later on. The following lemma, which is inspired by Lemma 2.1 in \cite{reiher-odd-cycles-locally-dense}, shows that such a conversion is possible in $\F_p^n$ (including for $n=1$) if the forms of the underlying system are distinct. We include a proof as to the best of our knowledge the precise version we shall need is not available in the literature.
\begin{lemma}[Conversion of maps to sets]\label{lem:maps-to-sets}
Let $G=\F_p^n$, $\conf\subseteq G^d$ a configuration given by a system of distinct forms $\Phi=(\phi_i)_{i\in [d]}$ in $r$ variables and $f\colon G\rightarrow [0,1]$. There exists a set $A\subseteq G$ with
\begin{align*}
t_\conf(A) + t_\conf(A^C) \leq t_\conf(f) + t_\conf(1-f) + \frac{\binom{d}{2}}{\vert G \vert}.
\end{align*}
\end{lemma}
\begin{proof}
We define the \textit{non-injective instances} of $\conf$ as
\begin{align*}
\niconf=\{v\in \conf\subseteq G^d: \exists i,j\in [d]\colon i\neq j \land v_i=v_j\}
\end{align*}
and refer to $\iconf=\conf\setminus \niconf$ as the \textit{injective instances} of $\conf$. Consider the map
\begin{align*}
\psi \colon \Map(G,[0,1]) \rightarrow \R \qquad g\mapsto \sum_{v\in \iconf} \prod_{i=1}^k g(v_i)+\sum_{v\in \iconf} \prod_{i=1}^k (1-g(v_i)).
\end{align*}
The map $\psi$ is continuous and defined on a compact space, therefore it attains its minimum. Among the functions $g$ that yield this minimum, let $g_0$ be such that the set $R(g)=\{x\in G:g(x)\notin \{0,1\}\}$ is of least cardinality. We show that $R:=R(g_0)$ is in fact empty. Assume to the contrary that there exists $a\in R$. Consider then for each $\eta\in \R$ the function
\begin{align*}
g_\eta\colon G\rightarrow \R \qquad x\mapsto \begin{cases}g_0(a)+\eta &\text{if } x=a,\\ g_0(x) &\text{otherwise.}\end{cases}
\end{align*}
Because $g_\eta(a)\notin \{0,1\}$, we have $\im (g_\eta)\subset [0,1]$ for sufficiently small $\eta$. It is easy to see that the function $\eta\mapsto \psi(g_\eta)$ is affine linear since $a$ can only appear once as a coordinate in every injective instance. But as $\psi(g_0)$ is minimal, the linear coefficient of $\eta\mapsto \psi(g_\eta)$ must in fact be 0. We may now choose $\eta$ such that $g_\eta(a)=0$ while $\psi(g_\eta)=\psi(g_0)$, contradicting the choice of $g_0$.\\\\
Now let $A=g_0\inv(1)$ so that $1_A=g_0$. We obtain
\begin{align*}
t_\conf(A)+t_\conf(A^C)&\leq\frac{1}{\vert \conf\vert}\lp \psi(1_A) + \vert \niconf\vert \rp \leq \frac{1}{\vert \conf\vert}\lp \psi(f) + \vert \niconf\vert \rp \leq t_\conf(f)+t_\conf(1-f) +  \frac{\vert \niconf \vert}{\vert \conf\vert}.
\end{align*}
Note that $\vert\conf\vert=\lv G\rv^r$, so that the proof is complete if for each fixed pair $\{i,j\}\in [d]^{(2)}$, the set $\{v\in \conf:v_i=v_j\}$ has size at most $\lv G\rv^{r-1}$. To see this, it suffices to apply \Cref{lemma:reparametrization} to $\{\phi_i,\phi_j\}$, if the two forms are linearly independent, or just to $\phi_i$, if $\phi_j$ is a multiple of $\phi_i$. 
\end{proof}
\subsection{Exponential sums and the Fourier transform}
All of our results are proved using discrete Fourier analysis and we shall fix here our notation and recall some elementary results. The reader familiar with these concepts is encouraged to skip the remainder of this section. For $p$, $x\in \Z_p$ and $\omega$ a $p$th root of unity, we write $\omega^x$ to mean $\omega$ taken to the power of any representative of the residue class $x$. If the exponent is given by some function $\vp$ of $x$, $\omega^{\vp(x)}$ is called a \textit{phase function} and $\vp(x)$ is its \textit{phase}.\\\\
More generally, for a finite Abelian group $G$, we will make use of its  \textit{characters}, i.e., group homomorphisms $\gamma \colon G\rightarrow \{z\in \C:\lv z\rv=1\}$, and the group $\hat{G}$ they form under pointwise multiplication. The neutral element of $\hat{G}$ is the constant function $1$ and the inverse of a character $\gamma$ is its pointwise complex conjugate $\ov{\gamma}$.\\\\
The reason that characters are useful for us is that they can detect the zero element of the group, as the following identity shows.
\begin{align}\label{eq:char-ortho-sum}
\sum_{\gamma\in \hat{G}} \gamma(x)=\begin{cases}\vert G\vert &\text{if } x=0,\\0&\text{otherwise.}\end{cases}
\end{align}
Given a function $f\colon G\rightarrow \C$ we define the \textit{Fourier transform} of $f$ as the function
\begin{align*}
\hat{f}\colon \hat{G} \rightarrow \C \qquad \gamma \mapsto \hat{f}(\gamma)=\E_{x\in G} \gamma(x)f(x).
\end{align*}
Note that for the trivial character 1 we have $\hat{f}(1)=\E(f)$. In a dual fashion, the constant function 1 has a particularly simple Fourier transform, namely,
\begin{align}\label{eq:fourier-transform-constant}
\hat{1}(\gamma)=\begin{cases}1&\text{if } \gamma = 1,\\0&\text{otherwise.}\end{cases}
\end{align}
If $f$ is real, and in fact only then, $\hat{f}(\gamma\inv)=\ov{\hat{f}(\gamma)}$ for every $\gamma\in \hat{G}$. Lastly, recall Plancherel's identity, which states that for all $f_1,f_2\in \Map(G,\C)$
\begin{align*}
\sum_{\gamma\in \hat{G}} \hat{f_1}(\gamma)\ov{\hat{f_2}(\gamma)}=\E_{x\in G} f_1(x)\ov{f_2(x)}.
\end{align*}
Functions that have large Fourier coefficients at frequencies other than $\gamma=1$ are thought of as being \textit{linearly structured} because they correlate with linear phase functions. For us, \textit{quadratically structured} functions will also be of interest, mainly due to the fact that they do not correlate with linear phase function. The following standard Gauss sum estimate, which we state without proof, formalizes this lack of correlation.
\begin{lemma}\label{lemma:quad-phase-vanish}
Let $G=\F_p^n$ for a prime $p$ and $n\in \N$. Suppose that $\omega$ is a $p$th root of unity and $a,c\in \F_p$, $b\in \F_p^n$ are such that $a$ and $b$ are not both $0$. Then 
\begin{align}\label{eq:quad-phase-vanish}
\lv\E_{x\in G} \omega^{ax^Tx+b^Tx+c} \rv \leq \vert G\vert^{-1/2}.
\end{align}
In particular,
\begin{align}\label{eq:mixed-phase-vanish}
\E_{y\in G}\lv\E_{x\in G} \omega^{ax^Tx+dy^Tx+b^Tx+c} \rv \leq \vert G\vert\inv + \vert G\vert^{-1/2}
\end{align}
if $(a,b,d)\neq (0,0,0)$.
\end{lemma}
\section{The muting method for constructing counterexamples}\label{sec:4ap-vector}
In this section, we introduce a method for proving that a given linear configuration is not common, which we call the \textit{muting method}. It can be understood as a generalization of Wolf's proof in \cite{wolf-4ap-uncommon} that the 4-AP itself is not common, which in turn is based on a construction by Gowers \cite{gowers-4ap-uniform-non-sidorenko}. The method is most easily explained in $G=\F_p^{n}$, where we construct functions $f\colon G\rightarrow [0,1]$ as follows. We take subspaces $V_1,V_2\subset\F_p^n$ that are spanned by distinct sets of standard basis vectors in $\F_p^n$ and two real-valued functions $f_1,f_2$ defined on $V_1$ and $V_2$, respectively. Writing $\pi_1$ and $\pi_2$ for the projections of $\F_p^n$ to $V_1$ and $V_2$, we set $f=1/2+ f_1\circ \pi_1\cdot f_2\circ \pi_2$. The arithmetic multiplicity of a configuration $\conf$ arising from a system of forms $\Phi=(\phi_1,\ldots,\phi_d)$, $\phi_i\colon G^r\rightarrow G$, is then given by
\begin{align}\label{eq:muting1}
t_\conf(f)&=\E_{v\in G^r} \prod_{i=1}^d \lp 1/2+ f_1\circ \pi_1(\phi_i(v))\cdot f_2\circ \pi_2(\phi_i(v))\rp\nonumber\\
&=2^{-d}+\sum_{k=1}^d 2^{k-d} \sum_{I\in [d]^{(k)}} \E_{v\in G^r} \prod_{i\in I}  f_1\circ \pi_1\circ \phi_i(v)\cdot f_2\circ \pi_2\circ \phi_i(v).
\end{align}
Because $f_1\circ \pi_1$ and $f_2\circ\pi_2$ only depend on the parts of $v$ in $V_1$ and  $V_2$ respectively, we can split the expectation according to these subspaces and separate $f_1$ from $f_2$ to obtain
\begin{align*}
t_\conf(f)=2^{-d}+\sum_{k=1}^d 2^{k-d} \sum_{I\in [d]^{(k)}} \lp \E_{v\in V_1^r} \prod_{i\in I}  f_1\circ \phi_i(v) \rp \cdot \lp \E_{v\in V_2^r} \prod_{i\in I}  f_2\circ \phi_i(v) \rp.
\end{align*}
The term $2^{-d}$ is the main term and together with the main term of $t_\conf(1-f)$ adds up to $2^{1-d}$, the threshold for a function to be common. The remaining sum adds up the {contributions of different subconfigurations corresponding to the subsets $\emptyset \neq I\subseteq [d]$. Note that both of the inner expectations are themselves arithmetic multiplicities, $t_{\Phi_I}(f_1)$ and $t_{\Phi_I}(f_2)$, for the system of forms $\Phi_I=(\phi_i)_{i\in I}$. The central idea of muting is to construct $f_1$ in such a way that $t_{\Phi_I}(f_1)$ is negative for a certain class of ``good'' subconfigurations of even size, and to choose $f_2$ such that $t_{\Phi_I}(f_2)$ is very small in modulus for precisely those subconfigurations that are not good. If we can achieve this, then $f$ disproves commonness.\\\\
In what follows, we will refer to $f_1$ and $f_1\circ \pi$ for some projection $\pi$ as the \textit{directional part of $f$}, while $f_2$ and $f_2\circ \pi$ will be called the \textit{muting part of $f$}.\\\\
In anticipation of the proof of \Cref{thm:ap4vector}, it is helpful to first analyze systems of forms containing two forms that are multiples of one another. The proof of the following result can be seen as a primitive application of the muting method, although the line between the directional and the muting part is rather blurry.
\begin{proposition}\label{prop:uncommon-multiple}
For each $d\in \N$ there exists an $\e>0$ such that the following holds. Let $G=\F_p^n$ for some prime $p$, let $n\in \N$ and let $\Phi\colon G^r\rightarrow G^d$ be a system of distinct forms (mod $p$) containing two forms $\phi_i, \phi_j$ such that $\phi_j=c\phi$ for some $c\in \F_p\setminus\{0,1\}$. Then there exists $f\colon G\rightarrow [0,1]$ such that
\begin{align}\label{eq:multiple-forms-uncommon}
t_\Phi(f)+t_\Phi(1-f)\leq 2^{1-d}-\e.
\end{align}
In particular, $\Phi$ is not fully common in $G$.
\end{proposition}
\begin{proof}
Let $r,d,p,n$ and $\Phi$ be given. We show that there exists $\e$ depending only on $d$ and satisfying \eqref{eq:multiple-forms-uncommon} by considering two different cases. Note that by the distinctness of the forms and $c\neq 0$ we have $p\neq 2$. For the first case, assume that there exist $j_1,j_2\in [d]$ such that $\phi_{j_1}=-\phi_{j_2}$. Then we take as the muting part a constant $\alpha$ and as the directional part the function
\begin{align*}
f_1\colon \F_p^n\rightarrow [-2,2] \qquad x \mapsto \lp i\omega^{1^Tx}-i\omega^{-1^Tx}\rp,
\end{align*}
where $i$ is the imaginary unit and $1$ is the vector in $\F_p^n$ with $1$ in every component. Note that $f_1$ takes only real values. Letting $f=1/2+\alpha\cdot f_1$ for some $\alpha\in [0,1/4)$ we obtain
\begin{align}\label{eq:ramsey-mult-multiple-form}
t_\Phi(f)+t_\Phi(1-f)&=2^{1-d}+\sum_{\substack{k=2\\k\text{ even}}}^d 2^{k+1-d}\alpha^k \sum_{I\in [d]^{(k)}} \E_{v\in G^r} \prod_{j\in I} f_1\circ \phi_j(v)\nonumber\\
&\leq 2^{1-d}+2^{3-d}\alpha^2 \sum_{I\in [d]^{(2)}} \E_{v\in G^r} \prod_{j\in I} f_1\circ \phi_j(v)+ 2^d\cdot 4^d\cdot \alpha^4.
\end{align}
It is now enough to show that the second term without the factor of $\alpha^2$ is less than $-C\alpha^2$ for some constant $C>0$ because we can then pick $\alpha$ small enough to make the last summand negligible. To achieve this, fix $I=\{j_1,j_2\}\subset [d]$ and observe that
\begin{align}\label{eq:lambda-formula-multiple-form}
\prod_{j\in I} f_1\circ \phi_j(v) &=\prod_{j\in I}\lp i\omega^{1^T\phi_j(v)}-i\omega^{-1^T\phi_j(v)}\rp\nonumber\\
&=-\sum_{\lambda\in \{\pm 1\}^2} (-1)^{(\lambda_1-\lambda_2)/2}\omega^{\lambda_1 1^T\phi_{j_1}(v)+\lambda_2 1^T\phi_{j_2}(v)}.
\end{align}
If $\phi_{i_2}$ is not a multiple of $\phi_{i_1}$, the two forms are linearly independent so we can apply \Cref{lemma:reparametrization} to $\phi_{i_1}$ and $\phi_{i_2}$ to obtain
\begin{align*}
\E_{v\in G^r} \prod_{j\in I} f_1\circ \phi_j(v)&=-\E_{v\in G^r}\sum_{\lambda\in \{\pm 1\}^2} (-1)^{(\lambda_1-\lambda_2)/2}\omega^{\lambda_1 1^Tv_1+\lambda_2 1^Tv_2}\\
&=-\sum_{\lambda\in \{\pm 1\}^2} (-1)^{(\lambda_1-\lambda_2)/2}\lp \E_{x\in G} \omega^{\lambda_1 1^Tv_1} \rp \lp \E_{y\in G}\omega^{\lambda_2 1^Tv_2}\rp,
\end{align*}
which is zero. If, on the other hand, $\phi_{i_1}=a\phi_{i_2}$ for $a\neq 1$, then we apply \Cref{lemma:reparametrization} just to $\phi_{i_1}$ to get
\begin{align*}
\E_{v\in G^r} \prod_{j\in I} f_1\circ \phi_j(v)&=-\E_{x\in G}\sum_{\lambda\in \{\pm 1\}^2} (-1)^{(\lambda_1-\lambda_2)/2}\omega^{\lambda_1 1^Tx+a\lambda_2 1^Tx}\\
&=-\sum_{\lambda\in \{\pm 1\}^2} (-1)^{(\lambda_1-\lambda_2)/2} \E_{x\in G}\omega^{(\lambda_1+a\lambda_2)1^Tx}.
\end{align*}
The inner average is only non-zero if $\lambda_1+a\lambda_2=0$, which is equivalent to $a=-\lambda_1/\lambda_2\in \{\pm 1\}$. Since $a\neq 1$, we then have that $a=-1$ and $\lambda_1=\lambda_2$. Therefore,
\begin{align*}
\E_{v\in G^r} \prod_{j\in I} f_1\circ \phi_j(v)&=\begin{cases}0&\text{if } a\neq -1,\\-2&\text{if } a=-1.
\end{cases}
\end{align*}
This means the contribution will not be positive for any pair $\{j_1,j_2\}$ and negative for at least the one pair where $\phi_{j_1}=-\phi_{j_2}$, the existence of which we had assumed. Substituting this into \eqref{eq:ramsey-mult-multiple-form} yields
\begin{align*}
t_\Phi(f)+t_\Phi(1-f)&\leq 2^{1-d}-2^{4-d}\alpha^2 + 8^d\cdot \alpha^4.
\end{align*}
Choosing $\alpha$ small enough allows us to find $\e$ satisfying \eqref{eq:multiple-forms-uncommon}.\\\\
We are now left with the case that there exists a pair $\{j_1,j_2\}$ with $\phi_{j_1}=c\phi_{j_2}$ for $c\in \F_p$, but no two forms are equal or additive inverses of one another. Fix $c$, let 
\begin{align*}
f_1\colon \F_p^n\rightarrow [-2,2] \qquad x \mapsto \lp \omega^{1^Tx}+\omega^{-1^Tx}-\omega^{c1^Tx}-\omega^{-c1^Tx}\rp
\end{align*}
and set $f=1/2+\alpha f_1$ as before. Once more, we only have to consider subconfigurations of size 2 and analogously to \eqref{eq:lambda-formula-multiple-form}, we obtain for a fixed pair $I=\{j_1,j_2\}$ that
\begin{align*}
\prod_{j\in I} f_1\circ \phi_j(v) &=\sum_{\lambda\in \{\pm 1,\pm c\}^2} \sigma(\lambda)\omega^{\lambda_1 1^T\phi_{j_1}(v)+\lambda_2 1^T\phi_{j_2}(v)},
\end{align*}
where $\sigma(\lambda)$ is $1$ if the components of $\lambda$ are both in $\{\pm 1\}$ or both in $\{\pm c\}$ and $-1$ if they contain one element from each of the two sets. If the forms $\phi_{j_1},\phi_{j_2}$ are independent, it follows from the same calculation as in the first case that the average over $v\in G^r$ vanishes. If $\phi_{j_2}=a\phi_{j_1}$ for $a\in \F_p\setminus \{0,\pm 1\}$, we reparametrize for $\phi_{j_1}$ and obtain
\begin{align*}
\E_{v\in G^r} \prod_{j\in I} f_1\circ \phi_j(v)&=\sum_{\lambda\in \{\pm 1, \pm c\}^2}\E_{x\in G} \sigma(\lambda)\omega^{(\lambda_1+a\lambda_2) 1^Tx}.
\end{align*}
If $\lambda_1/\lambda_2=-a\neq \pm 1$, then we know that the components of $\lambda$ cannot both lie in $\{\pm 1\}$ or both lie in $\{\pm c\}$, hence $\sigma(\lambda)=-1$. That means no subconfiguration of size 2 can make a positive contribution to the arithmetic multiplicity. At the same time, from the pair $\{j_1,j_2\}$ with $\phi_{j_2}=c\phi_{j_1}$ we obtain a negative contribution for $\lambda\in \{(c,-1),(-c,1)\}$. If we substitute this into \eqref{eq:ramsey-mult-multiple-form} and choose $\alpha$ sufficiently small, we find again $\e>0$ such that \eqref{eq:multiple-forms-uncommon} holds.
\end{proof}
Of course, the above proof could be adapted to yield that $\Phi$ is not common (instead of not fully common) by an application of \Cref{lem:maps-to-sets}, provided the group is sufficiently large.\\\\
We now go on to construct functions on $\F_p$ that have negative arithmetic multiplicity. These will be used as the directional part for the muting method as applied in the proofs of \Cref{thm:ap4vector} and \Cref{thm:ap4cyclic}.
\begin{lemma}\label{lemma:direction}
There exists $c\in \R$ such that for every prime $p\geq 5$, there exists a map $f\colon\F_p\rightarrow~[-1,1]$ with
\begin{align}\label{eq:ap4-neg-density}
\E_{x,y\in \F_p} f(x)f(x+y)f(x+2y)f(x+3y) \leq c.
\end{align}
Moreover, $f$ can be taken of the form $f=a1_{\pi([M]_0)}-b\cdot 1_{\pi(A)}$, where $\pi\colon \Z\rightarrow \F_p$ is the projection to residue classes, $M\in \N$, $a\in \{1,1/2\}$, $b\in \{2,3/2\}$ and $A\subset [M]_0$ is an arithmetic progression with common difference at most 5.
\end{lemma}
The construction of $f$ will depend on the prime $p$ and we will need to distinguish six different cases. We give the proof for the case $p\geq 200$ in full detail, but for the remaining ones we only provide the parameters $M,a,b,A$, which will allow the keen reader to verify the claim.\footnote{In an earlier arXiv version \cite{versteegen-4ap-uncommon-v2-arxiv} of this paper all six cases are calculated explicitly.} We point out that for $p=5$ the function $f$ is the same as in \cite{wolf-4ap-uncommon} and for larger $p$ we merely generalize the same construction.
\begin{proof}
There will be six cases to consider.\\\\
\textbf{Case 1:} $p=5$. We choose $M=4$, $a=1$, $b=2$ and $A=\{0\}$.\\\\
\textbf{Case 2:} $p=7$. We choose $M=6$, $a=1/2$, $b=3/2$ and $A=\{0\}$.\\\\
\textbf{Case 3:} $p=11$. We choose $M=6$, $a=1$, $b=2$ and $A=\{0,3\}$.\\\\
\textbf{Case 4:} $p=13$. We choose $M=7$, $a=1$, $b=2$ and $A=\{0,5\}$.\\\\
\textbf{Case 5:} $p\in [17,200]$. We choose $M=7$, $a=1$, $b=2$ and $A=\{0,5\}$.\\\\
\textbf{Case 6:} $p\geq 200$.
We choose $M=(p-1)/2>100$, $a=1$, $b=2$ and $A=\{5k:k\in \N_0, 5k\leq M\}$. Note that every arithmetic progressions in $\F_p$ that makes a non-zero contribution to the left hand side of \eqref{eq:ap4-neg-density} must arise from an integer 4-AP in $[M]_0$. For $m,n\in \Z$ with $0\leq m,m+3n\leq M$, we denote the contribution that the 4-AP $(m,m+n,m+2n,m+3n)$ makes by $\sigma(m,n)$. It is easy to see that
\begin{align}\label{eq:sigma-mod5}
\sigma(m,n)=\begin{cases}1&\text{if } n\equiv 0 \text{ or } n\equiv m \mod 5,\\
-1&\text{otherwise.}\end{cases}
\end{align}
We have
\begin{align*}
\sum_{x,y\in \F_p} f(x)f(x+y)f(x+2y)f(x+3y)&=\sum_{\substack{m,n\in \Z\\0\leq m,m+3n\leq M}} \sigma(m,n).
\end{align*}
For every pair $(m,n)$ where $n$ is negative, the 4-AP corresponding to this pair contains the same points as $(m-3n,-n)$ so we may free ourselves of all negative $n$ by bounding
\begin{align*}
\sum_{\substack{m,n\in \Z\\0\leq m,m+3n\leq M}} \sigma(m,n)\leq 2\sum_{\substack{m,n\in \N_0\\m+3n\leq M}} \sigma(m,n)=2\sum_{m=0}^M \sum_{n=0}^{\left\lfloor\frac{M-m}{3}\right\rfloor}\sigma(m,n).
\end{align*}
We write $r(m)$ for $\left\lfloor \left\lfloor\frac{M-m}{3} \right \rfloor/5\right\rfloor$ to split the sum over $n$ into blocks of length 5 and obtain
\begin{align*}
\sum_{m=0}^M \sum_{n=0}^{\left\lfloor\frac{M-m}{3}\right\rfloor}\sigma(m,n)&=\sum_{m=0}^M \lp \sum_{k=0}^{r(m)-1}\lp \sum_{i=0}^4\sigma(m,5k+i)\rp + \sum_{n=5r(m)}^{\left\lfloor\frac{M-m}{3}\right\rfloor}\sigma(m,n)\rp.
\end{align*}
Note that for each $m,k\in \N_0$ we have $\sum_{i=0}^4 \sigma(m,5k+i)=2-3=-1$ by \eqref{eq:sigma-mod5}. We can bound the above expression by
\begin{align*}
\sum_{m=0}^M \lp -r(m) + 2\rp\leq \sum_{m=0}^M \lp -\frac{M-m}{15} + 3\rp=3M-\frac{M^2}{15}+\frac{M(M+1)}{30}=\frac{91M-M^2}{30}.
\end{align*}
Since $M>100$, the expression above is negative. Concretely, we get
\begin{align*}
\E_{x,y\in \F_p} f(x)f(x+y)f(x+2y)f(x+3y)\leq \frac{2}{p^2}\cdot \frac{91M-M^2}{30} < \frac{1}{60p} - \frac{1}{120}.
\end{align*}
\end{proof}
A careful check of the six cases in the above proof shows that it is permissible to take $c=-2/199^2$ in the statement of \Cref{lemma:direction}. The important point for us, however, is that there exists a constant independent of $p$, which is crucial for the proof of \Cref{thm:ap4cyclic}.\\\\
We are now ready to prove \Cref{thm:ap4vector}, which we recall here for the convenience of the reader.
\uncommonapvector*
\begin{proof}
We may assume that no two forms are multiples of one another, else we would be done by \Cref{prop:uncommon-multiple}. For now, we let $n$ be any natural number. We shall impose the constraint on its size at the end of the argument. We let $f_1'$ be the function obtained in \Cref{lemma:direction} and precompose it with a projection $\pr_1\colon \F_p^{n+1}\rightarrow \F_p^n$ to the first coordinate to get $f_1=f_1'\circ\pr_1\colon \F_p^{n+1}\rightarrow [-1,1]$. Furthermore, we let
\begin{align*}
f_2'\colon \F_p^n \rightarrow [-4,4] \qquad x \mapsto \lp\beta\omega^{x^Tx+1^Tx} + \beta\omega^{-x^Tx-1^Tx} + \omega^{3x^Tx+3\cdot 1^Tx} + \omega^{-3x^Tx-3\cdot 1^Tx}\rp,
\end{align*}
where $\omega$ is a $p$th root of unity and $\beta\in (0,1]$ is a constant to be chosen later, which may depend on the matrix $M$ inducing $\Phi$, but not on $n$. Now set $f_2=f_2'\circ \pr_R$, where $\pr_R$ is the projection $\F_p^{n+1}\rightarrow \F_p^n$ onto the remaining $n$ coordinates, and finally let
\begin{align*}
f=\frac{1}{2}+\alpha f_1\cdot f_2.
\end{align*}
The constant $\alpha >0$ should always be small enough for $f$ to take only values in $[0,1]$ but will possibly be chosen even smaller to ensure $f$ satisfies \eqref{eq:goal-Fpn}.\\\\
Consider now the expanded arithmetic multiplicities
\begin{align*}
t_\Phi(f)+t_\Phi(1-f)&=\lp \frac{1}{2}\rp^{d-1}+\sum_{\substack{k=2\\k \text{ even}}}^d \lp\frac{1}{2}\rp^{d-k-1}\alpha^k \sum_{I\in [d]^{(k)}} \E_{v\in G^r} \lp\prod_{i\in I} f_1(\phi_i(v))f_2(\phi_i(v))\rp\\
&\leq \lp \frac{1}{2}\rp^{d-1}+\sum_{\substack{k=2\\k \text{ even}}}^4 \frac{\alpha^k}{2^{d-k-1}} \sum_{I\in [d]^{(k)}} \E_{v\in G^r} \lp\prod_{i\in I} f_1(\phi_i(v))f_2(\phi_i(v))\rp + \alpha^6 \cdot 8^d.
\end{align*}
Observe that for any $I\subseteq [d]$ we can always split the expectation over $v\in G^r$ as
\begin{align*}
\E_{v\in G^r} \lp\prod_{i\in I} f_1(\phi_i(v))f_2(\phi_i(v)) \rp=\lp \E_{u\in \F_p^r} \prod_{i\in I} f_1'(\phi_i(u)) \rp \lp \E_{w\in (\F_p^n)^r} \prod_{i\in I} f_2'(\phi_i(w)) \rp.
\end{align*}
By our choice of $f_1'$ we know that the first factor will be less than some negative constant $c$ whenever $I$ corresponds to a 4-AP. Therefore, it suffices to prove the following three statements to conclude our proof.
\begin{enumerate}
\item For $\vert I\vert =2$, we have
\begin{align}\label{eq:vector-size2-bound}
\lv \E_{w\in (\F_p^n)^r} \prod_{i\in I} f_2'(\phi_i(w)) \rv \leq \frac{16}{\sqrt{p^n}}.
\end{align}
\item If $\vert I\vert=4$ and $I$ does not correspond to a 4-AP, we have
\begin{align}\label{eq:vector-non-4ap-bound}
\lv \E_{w\in (\F_p^n)^r} \prod_{i\in I} f_2'(\phi_i(w)) \rv \leq \frac{244}{\sqrt{p^n}}+12\beta^3.
\end{align}
\item If $I$ corresponds to a 4-AP, we have
\begin{align}\label{eq:vector-4ap-bound}
\E_{w\in (\F_p^n)^r} \prod_{i\in I} f_2'(\phi_i(w)) \geq 2\beta^2-\frac{254}{\sqrt{p^n}}.
\end{align}
\end{enumerate}
Indeed, if the above inequalities hold, then there exist positive constants $C_0,C_1, C_2$ and $C_3$ depending only on $d$ such that
\begin{align}\label{eq:final-bound-4ap-vector}
t_\Phi(f)+t_\Phi(1-f)&\leq \lp \frac{1}{2}\rp^{d-1}+\sum_{\substack{k=2\\k \text{ even}}}^4 \frac{\alpha^k}{2^{d-k-1}} \sum_{I\in [d]^{(k)}} \E_{v\in G^r} \lp\prod_{i\in I} f_1(\phi_i(v))f_2(\phi_i(v))\rp + \alpha^6 \cdot 8^d\nonumber\\
&\leq \lp \frac{1}{2}\rp^{d-1}- C_0\alpha^4\beta^2 +\frac{C_1}{\sqrt{p^n}}+ C_2\alpha^4\beta^3  + C_3\alpha^6.
\end{align}
Choosing, in this order, $\beta$ sufficiently small, $\alpha$ sufficiently small and $n$ sufficiently large yields a positive $\e$ satisfying \eqref{eq:goal-Fpn}.\\\\
We first prove (i). For $I=\{i_1,i_2\}$ we have
\begin{align*}
\E_{w\in (\F_p^n)^r}\prod_{i\in I} f_2'(\phi_i(w))&=\\&\hspace*{-1.0cm}\E_{w\in (\F_p^n)^r}\sum_{\lambda\in \{\pm 1,\pm 3\}^2} \beta^{\rho_2(\lambda)}\omega^{\lambda_1(\phi_{i_1}(w)^T\phi_{i_1}(w)+1^T\phi_{i_1}(w))+\lambda_2(\phi_{i_2}(w)^T\phi_{i_2}(w)+1^T\phi_{i_2}(w))}
\end{align*}
where $\rho_2=\vert \{ j\in [2]:\lambda_j=\pm 1\} \vert$. Because $\phi_{i_1}$ and $\phi_{i_2}$ are independent we can reparametrize using \Cref{lemma:reparametrization} to obtain
\begin{align*}
\lv \E_{w\in (\F_p^n)^r}\prod_{i\in I} f_2'(\phi_i(w))\rv&=\lv\sum_{\lambda\in \{\pm 1,\pm 3\}^2}\E_{w\in (\F_p^n)^{r-2}}\E_{x,y\in \F_p^n} \beta^{\rho_2(\lambda)}\omega^{\lambda_1(x^Tx+1^Tx)+\lambda_2(y^Ty+1^Ty)}\rv\\
&\leq \sum_{\lambda\in \{\pm 1,\pm 3\}^2}\lv\E_{x\in \F_p^n} \beta^{\rho_2(\lambda)}\omega^{\lambda_1(x^Tx+1^Tx)}\rv,
\end{align*}
which is bounded by $\frac{16}{\sqrt{p^n}}$. This establishes \eqref{eq:vector-size2-bound}.\\\\
In order to prove (ii), let $I=\{i_1,i_2,i_3,i_4\}$. This time, the product $\prod_{i\in I} f_2'(\phi_i(w))$ will expand to  256 terms, one for each $\lambda\in \{\pm 1,\pm 3\}^4$. Each of these will be of the form
\begin{align}\label{eq:lambda-term-4ap}
\beta^{\rho_4(\lambda)}\omega^{\tau_\lambda(w)},
\end{align}
where $\rho_4(\lambda)=\vert \{ j\in [4]:\lambda_j=\pm 1\} \vert$ and
\begin{align*}
\tau_\lambda(w)&=\lambda_1(\phi_{i_1}(w)^T\phi_{i_1}(w)+1^T\phi_{i_1}(w))+\lambda_2(\phi_{i_2}(w)^T\phi_{i_2}(w)+1^T\phi_{i_2}(w))\\&\hspace*{2.0cm}+\lambda_3(\phi_{i_3}(w)^T\phi_{i_3}(w)+1^T\phi_{i_3})+\lambda_4(\phi_{i_4}(w)^T\phi_{i_4}(w)+1^T\phi_{i_4}(w)).
\end{align*}
Denote by $s$ the dimension of the linear span of $\phi_{i_1},\ldots,\phi_{i_4}$ in the space of linear forms from $(\F_p^n)^r$ to $\F_p^n$. By the assumption of pairwise linear independence, $s$ is at least 2. If $s=4$, the forms are linearly independent overall and the linear coefficient of the phase will not vanish. If we reparametrize for $\phi_{i_1}$, we can apply \Cref{lemma:quad-phase-vanish} to bound the average of \eqref{eq:lambda-term-4ap} over $w_1\in G$ by $p^{-n/2}$ to get \eqref{eq:vector-non-4ap-bound}.\\\\
If $s=3$, we may assume without loss of generality that $\phi_{i_1},\phi_{i_2},\phi_{i_3}$ are linearly independent. We apply \Cref{lemma:reparametrization} to get $\phi_{i_1}',\ldots,\phi_{i_4}'$ with $\phi_{i_1}'(w)=w_1,\ldots,\phi_{i_3}'(w)=w_3$. By pairwise linear independence, we know that $\phi_{i_4}'$ depends on at least two of the variables $w_1,w_2,w_3$. Assume it depends on $w_1$ and $w_2$. Therefore, the product $w_1^Tw_2$ will necessarily appear in $(\phi_{i_4}')^T\phi_{i_4}'(w)$ with a non-zero coefficient. At the same time, $w_1^Tw_2$ does not appear in any of the other squares, which are just $w_1^Tw_1,w_2^Tw_2$ and $w_3^Tw_3$. By \eqref{eq:mixed-phase-vanish} from \Cref{lemma:quad-phase-vanish}, averaging over $w_1,w_2$ then yields \eqref{eq:vector-non-4ap-bound}.\\\\
This leaves us with the case $s=2$. By reordering $i_1,\ldots,i_4$ and redefining $\tau_\lambda$ accordingly, we can achieve the following. If there exists $\lambda\in \{\pm 1\}^4$ such that $\tau_\lambda=0$ and exactly two components are $1$, then we also have $\tau_{(-1,1,1,-1)}=0$. If such a $\lambda$ does not exist, but instead there exists $\lambda\in \{\pm 1,\pm 3\}^4$ such that all coordinates are distinct and $\tau_\lambda=0$, then $\tau_{(-3,3,1,-1)}=0$ as well.\footnote{Note that we cannot insist that $\lambda=(-1,1,1,-1)$ or $\lambda=(-3,3,1,-1)$ respectively. There may well be several $\lambda$ such that $\tau_\lambda=0$.} Furthermore, we may reparametrize the forms using \Cref{lemma:reparametrization} such that $\phi'_{i_1}(w)=w_1=:x$ and $\phi'_{i_2}(w)=w_2=:y$. The remaining two forms can then be expressed as $\phi'_{i_3}(w)=a_1x+b_1y$ and $\phi'_{i_4}(w)=a_2x+b_2y$, where none of $a_1,b_1,a_2,b_2$ can be 0 and $(a_1,b_1)$ and $(a_2,b_2)$ are not multiples of one another.\\\\ For each $\lambda\in \{\pm 1,\pm 3\}$, we let
\begin{align*}
\tau'_{\lambda}(x,y)&=(\lambda_1+a_1^2\lambda_3+a_2^2\lambda_4) x^Tx+(\lambda_2+b_1^2\lambda_3+b_2^2\lambda_4)(y^Ty)+2(a_1b_1\lambda_3+a_2b_2\lambda_4)x^Ty\\
&\hspace*{2.0cm}+(\lambda_1+a_1\lambda_3+a_2\lambda_4)1^Tx + (\lambda_2+b_1\lambda_3+b_2\lambda_4)1^Ty,
\end{align*}
so that
\begin{align*}
\E_{w\in (\F_p^n)^r} \beta^{\rho_4(\lambda)}\omega^{\tau_\lambda(w)} = \E_{x,y} \beta^{\rho_4(\lambda)}\omega^{\tau'_\lambda(x,y)}.
\end{align*}
Whenever at least one of the coeffiecients of $x^Tx,y^Ty,x^Ty,x$ or $y$ in $\tau'_\lambda(x,y)$ is non-zero, we can apply \Cref{lemma:quad-phase-vanish} to bound the above average over $x,y$ by $p^{-n/2}$. We shall now determine when this is the case. The condition that all coefficients are 0 is expressed by the following homogeneous system of equations over $\F_p$
\begin{align}\label{eq:big-system}
\begin{pmatrix}
1&0&a_1^2&a_2^2\\0&1&b_1^2&b_2^2\\0&0&2a_1b_1&2a_2b_2\\1&0&a_1&a_2\\0&1&b_1&b_2
\end{pmatrix}
\begin{pmatrix}
\lambda_1\\\lambda_2\\\lambda_3\\\lambda_4
\end{pmatrix}=\begin{pmatrix}
0\\0\\0\\0\\0
\end{pmatrix},
\end{align}
together with the additional constraint that $\lambda\in \{\pm 1,\pm 3\}^4$. We transform the above matrix into
\begin{align*}
\begin{pmatrix}
1&0&a_1^2&a_2^2\\0&1&b_1^2&b_2^2\\0&0&2a_1b_1&2a_2b_2\\0&0&a_1-a_1^2&a_2-a_2^2\\0&0&b_1-b_1^2&b_2-b_2^2
\end{pmatrix}.
\end{align*}
This system has non-trivial solutions if and only if the third, fourth and fifth row are each multiples of one another (for this analysis we consider all matrix entries as elements of $\F_p$). Assuming that this is the case, we can show that  $a_2-a_2^2\neq 0$. Indeed, suppose to the contrary that $a_2-a_2^2\neq 0$. Then $a_1-a_1^2=0$ as well because otherwise the system would have full rank (as $2a_2b_2\neq 0$). Because $a_1,a_2\neq 0$, we obtain $a_1=a_2=1$ and the system reduces to
\begin{align*}
\begin{pmatrix}
1&0&1&1\\0&1&b_1^2&b_2^2\\0&0&2b_1&2b_2\\0&0&b_1-b_1^2&b_2-b_2^2
\end{pmatrix}.
\end{align*}
For this system not to have full rank we would need to have
\begin{align*}
2b_1(b_2-b_2^2)&=2b_2(b_1-b_1^2),
\end{align*}
implying that $b_1=b_2$ and hence $(a_1,b_1)=(a_2,b_2)$, which cannot be the case. By the same argument, we can conclude that $b_2-b_2^2\neq 0$.\\\\
We can extract more information from the fact that the system does not have full rank. Notice that we must have
\begin{align*}
\frac{a_1-a_1^2}{a_2-a_2^2}=\frac{a_1b_1}{a_2b_2}=\frac{b_1-b_1^2}{b_2-b_2^2},
\end{align*}
from which it follows that
\begin{align}
\frac{1-a_1}{1-a_2}&=\frac{b_1}{b_2}\Rightarrow (1-a_1)b_2=(1-a_2)b_1 \Rightarrow b_2-a_1b_2=b_1-a_2b_1\label{eq:diff-one}\\
\frac{1-b_1}{1-b_2}&=\frac{a_1}{a_2}\Rightarrow (1-b_1)a_2=(1-b_2)a_1 \Rightarrow a_1-b_2a_1=a_2-b_1a_2.\label{eq:diff-two}
\end{align}
After subtracting \eqref{eq:diff-one} from \eqref{eq:diff-two}, we obtain that $b_2-b_1=a_1-a_2$.
Note that the difference $c:=a_2-a_1=b_1-b_2$ cannot be 0 because $(a_1,b_1)\neq (a_2,b_2)$. Inserting $a_1=a_2+c$ and $b_2=b_1+c$ into the last identity of \eqref{eq:diff-one} yields
\begin{align*}
b_1+c-(a_2+c)(b_1+c)=b_1-a_2b_1\Rightarrow c(a_2+b_1+c-1)=0\Rightarrow b_1=1-a_2-c.
\end{align*}
Renaming $a:=a_2$, we have $a_1=a+c$, $b_1=1-a-c$ and $b_2=1-a$. Substituting this into \eqref{eq:big-system} and dividing the third equation by 2, we arrive at the equation
\begin{align}\label{eq:big-system2}
\begin{pmatrix}
1&0&(a+c)^2&a^2\\0&1&(1-a-c)^2&(1-a)^2\\0&0&(a+c)(1-a-c)&a(1-a)\\1&0&a+c&a\\0&1&1-a-c&1-a
\end{pmatrix}
\begin{pmatrix}
\lambda_1\\\lambda_2\\\lambda_3\\\lambda_4
\end{pmatrix}=\begin{pmatrix}
0\\0\\0\\0\\0
\end{pmatrix}.
\end{align}
Let us split $\{\pm 1, \pm 3\}^4$ into four classes, namely
\begin{align*}
A_1=\{&\lambda\in \{\pm 1,\pm 3\}^4:\lambda_1+\lambda_2+\lambda_3+\lambda_4\neq 0\},\\
A_2=\{&(-3,1,1,1),(1,-3,1,1),(1,1,-3,1),(1,1,1,-3),(3,-1,-1,-1),(-1,3,-1,-1),\\
&(-1,-1,3,-1),(-1,-1,-1,3)\},\\
A_3=\{&(3,3,-3,-3), (-3,-3,3,3),(3,-3,3,-3),(3,-3,-3,3),(-3,3,3,-3),(-3,3,-3,3),\\
&(1,1,-1,-1), (-1,-1,1,1),(1,-1,1,-1),(-1,1,-1,1),(1,-1,-1,1),(-1,1,1,-1)\},\\
A_4=\{&\lambda\in \{\pm 1,\pm 3\}^4:\{\lambda_1,\lambda_2,\lambda_3,\lambda_4\}=\{\pm 1,\pm 3\}\}.
\end{align*}
We bound the left-hand side of \eqref{eq:vector-non-4ap-bound} by
\begin{align*}
\lv \E_{w\in (\F_p^n)^r} \prod_{i\in I} f_2'(\phi_i(w)) \rv \leq \sum_{\lambda\in \{\pm 1,\pm 3\}^4} \lv \E_{x,y\in F_p^n} \beta^{\rho_4(\lambda)} \omega^{\tau'_\lambda(x,y)} \rv.
\end{align*}
Adding the last two equations in \eqref{eq:big-system2} gives $\lambda_1+\lambda_2+\lambda_3+\lambda_4=0$, meaning that for each $\lambda\in A_1$ the phase $\tau'_\lambda$ is non-trivial and hence the average can be bounded by $p^{-n/2}$. For $\lambda \in A_2$ on the other hand, we always have $\rho_4(\lambda)=3$ so that we can bound the average by $\beta^3$.\\\\
Consider now $\lambda\in A_3$. If $\lambda$ were a solution to \eqref{eq:big-system2} for some $a,c$, we could assume without loss of generality, by linearity and homogenity of \eqref{eq:big-system2} as well as our earlier reordering of $i_1,\ldots,i_4$, that $\lambda=(-1,1,1,-1)$. We could then infer from the fourth equation that $-1+a+c-a=0$ and hence $c=1$. Substituting this into the third equation would give
\begin{align*}
-(a+1)a-a(1-a)=-2a=0,
\end{align*}
implying, as $p\neq 2$, that $a=0$, which is forbidden.\\\\
This leaves $\lambda \in A_4$. Because we have already established that \eqref{eq:big-system2} has no solution in $A_3$, we know by our earlier reordering that if $\tau'_\lambda=0$, then $\tau'_{(-3,3,1,-1)}=0$. Let us assume for simplicity that $\lambda=(-3,3,1,-1)$. Let $a,c\in \F_p^*$ be such that $\lambda$ is a solution to \eqref{eq:big-system2}. By the fourth equation, $-3+a+c-a=0$, which implies $c=3$. Substituting this into the first equation gives $-3+(a+3)^2-a^2=0$ and thus $a=-1$.\\\\
By definition of $a$ and $c$, this means that $a_1=2,b_1=-1,a_2=-1,b_2=2$. Therefore, the configuration in question is $(x,y,2x-y,-x+2y)$ which is a 4-AP with common difference $x-y$. But this is forbidden by our assumption in (ii).\\\\
We have thus bounded $\lv \E_{x,y\in F_p^n} \beta^{\rho_4(\lambda)} \omega^{\tau'_\lambda(x,y)} \rv$ by $\beta^3$ for each of the twelve $\lambda$ from $A_2$ and by $1/\sqrt{p^n}$ for the 244 other $\lambda$ from $A_1$, $A_3$ and $A_4$. This yields \eqref{eq:vector-non-4ap-bound}.\\\\
Lastly, let us show \eqref{eq:vector-4ap-bound} for $I=\{i_1,i_2,i_3,i_4\}$ such that the forms $(\phi_{i_1},\ldots,\phi_{i_4})$ form a 4-AP. By reordering and reparametrizing, we can achieve  that for $w=(x,y,w_3,\ldots)\in G^r$ we have $\phi_1(w)=x$, $\phi_2(w)=x+y$, $\phi_3(w)=x+2y$ and $\phi_4(w)=x+3y$. Again, we expand
\begin{align*}
\E_{w\in (\F_p^n)^r} \prod_{i\in I} &f_2'(\phi_i(w)) \\&= \sum_{\lambda\in \{\pm 1,\pm 3\}^4} \E_{x,y\in F_p^n} \beta^{\rho_4(\lambda)} \omega^{\lambda_1x^Tx+\lambda_2(x+y)^T(x+y)+\lambda_3(x+2y)^T(x+2y)+\lambda_4(x+3y)^T(x+3y)}.
\end{align*}
For those $\lambda$ where the phase does not vanish, the expectation over $x,y$ will be bounded in modulus by $1/\sqrt{p^n}$. For those $\lambda$ where the phase does vanish, we are left with the expectation of the constant function 1. This happens for $(\lambda_1,\lambda_2,\lambda_3,\lambda_4)=\pm(1,-3,3,-1)$, for which $\rho_4(\lambda)=2$. Taking these two cases together, we have
\begin{align*}
\E_{w\in (\F_p^n)^r} &\prod_{i\in I} f_2'(\phi_i(w))\geq 2\beta^2 -\frac{254}{\sqrt{p^n}}.
\end{align*}
This is \eqref{eq:vector-4ap-bound}, concluding the proof.
\end{proof}
One might wonder whether the introduction of the factor $\beta$ in the muting part $f_2$ in the proof above is really necessary, that is, whether it is possible to show that no $\lambda$ from $A_2$ could ever be a solution to \eqref{eq:big-system2}. The answer to this question is negative. For example, $(1,-3,1,1)\in A_2$ is a solution to \eqref{eq:big-system2} when $p=7$, $a=2$ and $c=2$. But the subconfiguration described by these parameters is not a 4-AP either. Indeed, the tuple $(0,1,-1,-3)$ is an instance of the subconfiguration but does not form a 4-AP in $\F_7$, not even after reordering.\\\\
\Cref{cor:ap4vector} follows from \Cref{thm:ap4vector} by an application of \Cref{lem:maps-to-sets}. For this argument it is important that the constant $\e$ in \Cref{thm:ap4vector} be independent of $n$. For any given configuration, we could extract an explicit lower bound for $\e$ from \eqref{eq:final-bound-4ap-vector}, but this is a somewhat pointless exercise as for this given configuration there is no reason to believe that our construction should come anywhere close to what can be achieved by a more careful application of the muting method.
\section{The muting method in cyclic groups}\label{sec:4ap-cyclic}
We now turn to the case of cyclic groups. Working in $\F_p^n$ was convenient because we were able to let the directional and muting part depend on orthogonal subspaces so that they would not interfere with each other. The group $\Z_p$ does not have even subgroups, let alone orthogonal subspaces, so we cannot hope to reproduce this approach directly. Let us instead describe in more abstract terms why the two parts did not interfere with each other and see how this might transfer to a cyclic group. The key feature of our construction was that the directional part and the muting part (or rather their products over the different forms) are not strongly correlated with one another except for our ``good'' subconfigurations. Concretely, if the first coordinate of $\F_p^{n+1}$ was fixed, the directional part was constant while the muting part was a polynomial phase function. In keeping with this abstract viewpoint, we could say that \Cref{lemma:quad-phase-vanish} states that constant functions and polynomial phase functions are only very weakly correlated as long as the phase polynomial is non-trivial. Our strategy for the cyclic group will be to modify the muting part in such a way that it will be uncorrelated not only with constant functions but also with the indicator functions of the arithmetic progressions arising from \Cref{lemma:direction}, which form our directional part.\\\\ 
Arithmetic progressions are linear in nature. In particular, the indicator function of an arithmetic progression is well approximated by a function with few non-zero Fourier coefficients, and such a function can be written as a sum of a small number of linear phase functions. This is helpful, because yet another way to read \Cref{lemma:quad-phase-vanish} is that properly quadratic phase functions are not correlated with linear phase functions. We might try to make use of this by showing that the phase we called $\tau_\lambda$ in the proof of \Cref{thm:ap4vector} is in fact properly quadratic and not linear. However, this would involve letting go of the last two equations from \eqref{eq:big-system}, and it seems impossible to complete the proof without them. So instead of trying to remove the linear phase in the muting part, we are going to take a closer look at the linear structure of arithmetic progressions. It turns out that the directional part can only correlate strongly with linear phase functions for which the linear coefficient is a $C$-fraction. Thus, we can construct the muting part to have a linear coefficient that is not a $C$-fraction to avoid correlation.\\\\
The following lemma should be thought of as a generalization of \Cref{lemma:quad-phase-vanish}. In its proof, we shall distinguish carefully between elements of $\Z$ and their residue classes in $\Z_p$ but allow ourselves one abuse of notation for the sake of readability: if $g\in \Z_p$ and $k\in \Z$, then $g+k$ is the sum $g$ and the residue class of $k$ mod $p$ and $gk$ is their product.
\begin{lemma}\label{lem:phase-vanish-extended}
Let $p$ be a prime, $a,b,c,d\in \Z$ and assume we have for each $y\in \Z_p$ an arithmetic progression $A_y=\{t_y+js_y:0\leq j <L_y\}$, where $t_y\in \Z_p$ is the starting point of the progression $A_y$, $s_y\in \Z_p\setminus \{0\}$ is its common difference and $L_y\leq p$ is its length. Furthermore, let $C$ be a constant satisfying $4C^4<p$ and assume that for each $y\in \Z_p$ the common difference $s_y$ is a $C$-fraction and coprime to $p$. If $d$ is coprime to $p$, we have 
\begin{align}\label{eq:ext-vanishing-phase-d}
\E_{y\in \Z_p}\lv \E_{x\in \Z_p} \omega^{ax^2+bx+c+dxy} 1_{A_y}(x)\rv \leq 8C^2\frac{\log (p)+1}{p}.  
\end{align}
If $d$ is instead a multiple of $p$, we define\footnote{If we suppress $y$ by averaging only over $x$, then it does not matter which $A_y$ we consider. We simply choose $A_0$ for concreteness.} $A=A_0,s=s_0,t=t_0$ and $L=L_0$. If $a$ is coprime to $p$, we have
\begin{align}\label{eq:ext-vanishing-phase-a}
\lv \E_{x\in \Z_p} \omega^{ax^2+bx+c} 1_A(x)\rv \leq 2\frac{\log (p)+1}{\sqrt{p}}.
\end{align}
Finally, assume that $a$ and $d$ are multiples of $p$ but $b$ is not. Assume further that $b$ is of the form $b=b'q$, where $q=\lfloor \sqrt{p}\rfloor$ and $b'$ is a $C$-fraction. Then we have
\begin{align}\label{eq:ext-vanishing-phase-b}
\lv \E_{x\in \Z_p} \omega^{bx+c} 1_A(x)\rv \leq 2\frac{C^2+1}{\sqrt{p}}.
\end{align}
\end{lemma}
\begin{proof}
Fix a $y\in \Z_p$ and let $f_y(x)=\omega^{ax^2+bqx+c+dxy}$. We denote by $\g\colon \Z_p\rightarrow \hat{\Z_p}$ the canonical isomorphism from $\Z_p$ to $\hat{\Z_p}$. That is, if $g,h\in \Z_p$, then $\g(g)(h)=\omega^{gh}$. By Plancherel, we have that
\begin{align*}
\E_{x\in \Z_p} \omega^{ax^2+bqx+c+dxy} 1_{A_y}(x) = \sum_{\gamma\in \hat{Z_p}} \hat{f_y}(\gamma)\ov{\hat{1_{A_y}}(\gamma)}=\sum_{g\in \Z_p} \hat{f_y}(\g(g))\ov{\hat{1_{A_y}}(\g(g))}.
\end{align*}
To bound this expression, we investigate the Fouriers coefficients of $1_{A_y}$. We already know that $\hat{1_{A_y}}(\g(0))=\E(1_{A_y})=L_y/p\leq 1$. For a residue class $g\in \Z_p$, let $r(g)$ denote the representative of $g$ with the least modulus. For $g\neq 0$, we have
\begin{align*}
\lv\hat{1_{A_y}}(\g(g))\rv&=\frac{1}{p}\lv \sum_{j=0}^{L_y-1} \omega^{g(t_y+js_y)}\rv=\frac{1}{p}\lv \frac{1-\omega^{L_yg s_y}}{1-\omega^{gs_y}}\rv\leq \frac{2}{p}\lv 1-\exp\lp \frac{2\pi i}{p} r(g s_y)\rp \rv\inv.
\end{align*}
To bound this further from above, we would like to bound the distance of $\exp\lp 2\pi ir(g s_y)/p\rp$ from 1. We set $m_{g}:=r(gs_y)$. If $\lv m_{g}\rv>p/4$, then the argument $2\pi ir(gs_y)/p$ is at least $\pi/2$ away from a multiple of $2\pi$. In that case, the real part is negative, giving $\lv \hat{1_{A_y}}(\g(g))\rv\leq 2/p$. If the argument is within $\pi/2$ of a multiple of $2\pi$, i.e., $\lv m_{g}\rv< p/4$, we use the estimate $\sin(\alpha)\geq \alpha/2$ to bound
\begin{align*}
\lv 1 - \exp\lp \frac{2\pi i}{p} r(g s_y)\rp \rv \geq \lv \sin\lp \frac{2\pi}{p}m_{g}\rp\rv\geq \frac{\pi \vert m_{g}\vert}{p}.
\end{align*}
Hence if $g\neq 0$ (and therefore $m_{g}\neq 0$), then $\lv \hat{1_{A_y}}(\g(g))\rv \leq 1/\vert m_{g}\vert$.\\\\
What can we say about $\hat{f_y}$? \Cref{lemma:quad-phase-vanish} tells us that if $a$ is coprime to $p$, then all Fourier coefficients of $f_y$ will be bounded in modulus by $1/\sqrt{p}$. Because the map $g\mapsto m_{g}$ is injective, we have
\begin{align*}
\lv \sum_{g\in \Z_p} \hat{f_y}(g)\ov{\hat{1_{A_y}}(g)} \rv \leq \frac{1}{\sqrt{p}}+ \sum_{\substack{g\in \Z_p\setminus\{0\}\\ \vert m_{g} \vert < p/4}} \frac{1}{\sqrt{p}\cdot \vert m_{g}\vert} + \sum_{\substack{g\in \Z_p\\ \vert m_{g} \vert > p/4}} \frac{2}{\sqrt{p}^3}< 2\frac{\log(p)+1}{\sqrt{p}}.
\end{align*}
Letting $y=0$ gives \eqref{eq:ext-vanishing-phase-a}, while taking the average over $y\in \Z_p$ gives \eqref{eq:ext-vanishing-phase-d} under the assumption that $a$ is coprime to $p$. If, on the other hand, $a\equiv 0$ mod $p$, we have
\begin{align*}
\hat{f_y}(\g(g)) = \begin{cases}\omega^c & \text{if } g=-(b+dy),\\
0&\text{otherwise.}\end{cases}
\end{align*}
We set $g_y:=-(b+dy)$ and $m_y:=r(g_ys_y)$. By our earlier analysis of $\hat{1_{A_y}}$, we have 
\begin{align}\label{eq:ext-vanishing-d-inter}
\E_{y\in \Z_p}\lv \E_{x\in \Z_p} f_y(x)1_{A_y}(x) \rv&= \E_{y\in \Z_p} \lv \hat{1_{A_y}}(\g(g_y))\rv\nonumber\\
&\leq\frac{1}{p}\lp 1 + \sum_{\substack{y\in \Z_p\\ 0<\lv m_y\rv < p/4}} \frac{1}{\lv m_y\rv} + \sum_{\substack{y\in \Z_p\\ p/4<\lv m_y\rv < p/2}} \frac{2}{p} \rp.
\end{align}
If $d$ is coprime to $p$, then $y\mapsto g_y$ is a bijection and therefore no two elements $y\neq y'$ with $s_y=s_{y'}$ satisfy $m_y=m_{y'}$. But since each common difference $s_y$ is a $C$-fraction, by \Cref{lem:C-fraction-properties} there are at most $4C^2$ distinct common differences. This means that for each integer $m$ we have $\lv \{y\in\Z_p:m_y=m\}\rv\leq 4C^2$. The right-hand side of \eqref{eq:ext-vanishing-d-inter} is therefore bounded above by
\begin{align*}
\frac{1}{p}\lp 1+\sum_{\substack{m\in \Z\\ \vert m\vert<p/4}} \frac{4C^2}{\vert m\vert} + \sum_{\substack{m\in \Z\\ m/4< \vert g\vert<p/2}} \frac{8C^2}{p}\rp\leq 8C^2\frac{\log(p)+1}{p},
\end{align*}
whence \eqref{eq:ext-vanishing-phase-d} follows.
It remains to show \eqref{eq:ext-vanishing-phase-b} under the assumption that $a,d\equiv 0$ mod $p$ and $b=b'q$ is coprime to $p$. With $f=f_0$ we have
\begin{align*}
\lv \E_{x\in \Z_p} \omega^{bx+c} 1_A(x)\rv = \lv \hat{1_A}(-\g(b))\rv \leq \frac{2}{\vert r(-bs)\vert}=\frac{2}{\vert r(b'qs)\vert}.
\end{align*}
So all that remains to be done is to show that $b'qs$ cannot be too close to a multiple of $p$. Let $y_b,y_s\in \Z$ and let $z_b,z_s\in \N$ be bounded in modulus by $C$ such that $b'z_b\equiv y_b$ and $r(s)z_s\equiv y_s$ mod $p$. We have
\begin{align*}
r(b'qs)\equiv b'qr(s) \Leftrightarrow r(b'qs)z_bz_s\equiv y_by_sq \qquad \mod p.
\end{align*}
By assumption, $C^2$ is less than $\sqrt{p}/2$. If $\lv r(b'qs)\rv<\sqrt{p}$, both sides of the above equation mod $p$ are smaller than $p/2$ so that we must have $r(b'qs)z_bz_s=y_by_sq$. Either way, we have $\lv r(b'qs)\rv\geq q/C^2$ and hence
\begin{align*}
\lv \E_{x\in \Z_p} \omega^{bx+c} 1_A(x)\rv = \lv \hat{1_A}(-\g(b))\rv\leq \frac{2C^2}{q}\leq \frac{2C^2}{\sqrt{p}-1}\leq \frac{2C^2+1}{\sqrt{p}}.
\end{align*}
\end{proof}
In our proof of \Cref{thm:ap4cyclic} we will not distinguish between the three bounds \eqref{eq:ext-vanishing-phase-d}, \eqref{eq:ext-vanishing-phase-a} and \eqref{eq:ext-vanishing-phase-b} and simply observe that they are all smaller than
\begin{align*}
8C^2\frac{\log(p)+1}{\sqrt{p}}
\end{align*}
which approaches 0 as $p$ tends to infinity.\\\\
One crucial piece of preparation for the proof of \Cref{thm:ap4vector} is still missing. After all, it is not enough to see that the muting part does not correlate with the plain indicator functions. Rather, we need it to not correlate with a product of indicator functions precomposed by reparametrized linear forms. We will handle this by singling out one coordinate of $G^r$, call it $x$, that we average over. Fixing the other variables, we could say that we are now averaging in $x$ over the muting part multiplied by a product of indicator functions of scaled and translated arithmetic progressions. A product of indicator functions of sets is the indicator function of the intersection of the sets. We would like to apply \Cref{lem:phase-vanish-extended} at this point, but unfortunately, intersections of APs in $\Z_p$ are not typically APs again. This is why we need the following lemma, which allows us to split an AP in $\Z_p$ whose common difference is a $C$-fraction into a bounded number of projections of integer APs in the interval $[p-1]_0$ with a shared common difference of size at most $C$. These special APs have the benefit that any intersection of them yields again an AP.
\begin{lemma}[Splitting lemma]\label{lem:splitting}
Let $p$ be a prime and $A=\{t+sj\in \Z:0<j<L\}$ an arithmetic progression in $\Z$ whose common difference $s$ is a $C$-fraction coprime to $p$. Then there exist an integer $m\leq 3C$ and disjoint integer arithmetic progressions $A_1,\ldots,A_m\subseteq [0,p-1]$, all with the same  common difference $0<s'<C$, such that 
\begin{align*}
\pi\lp\bigcup_{i=1}^m A_i\rp =\pi (A),
\end{align*}
where $\pi\colon \Z\rightarrow \Z_p$ is the projection onto residue classes mod $p$.
\end{lemma}
\begin{proof}
Since $s$ is a $C$-fraction, there are integers $y$ and $z>0$ bounded by $C$ such that $sz\equiv y$ mod $p$. We assume without loss of generality that $s$ is positive and hence that $y$ is, too. We observe that
\begin{align*}
\pi(A)=\bigcup_{i=0}^{z-1} \pi\lp \{ (t+is)+jy: j\in \Z\colon 0<i+jz<L\}\rp.
\end{align*}
That is, we can write $\pi(A)$ as the projection of a union of $z$ arithmetic progressions with common difference $y$, each of which is contained in an interval of length at most $yL/z< yp/z$. We can split each of these progressions into $\lfloor y/z\rfloor +2$ segments, each of which is contained in $[kp,(k+1)p-1]$ for some $k\in \Z$. This gives a total of at most $(y/z+2)z=y+2z\leq 3C$ progressions that can be translated by a multiple of $p$ into $[0,p-1]$, proving the claim.
\end{proof}
We are now ready to prove \Cref{thm:ap4cyclic}.
\uncommonapcyclic*
\begin{proof}
Let $p>200$ be a prime. As before we may assume that no two forms are multiples of one another. We choose once more $f_1\colon \Z_p\rightarrow [-1,1]$ as given by \Cref{lemma:direction} for the case $p> 200$ as the directional part.\\\\
We also take essentially the same muting part as in the proof of \Cref{thm:ap4vector}, but scale the linear part of the phase by $q=\lfloor \sqrt{p}\rfloor$. That is,
\begin{align*}
f_2\colon \Z_p \rightarrow [-4,4] \qquad x \mapsto \lp\beta\omega^{x^2+qx} + \beta\omega^{-x^2-qx} + \omega^{3x^2+3qx} + \omega^{-3x^2-3qx}\rp,
\end{align*}
where $\beta\leq 1$ will be determined later. As before, we set
\begin{align*}
f=\frac{1}{2}+\alpha f_1\cdot f_2,
\end{align*}
and thereby obtain once more
\begin{align*}
t_\Phi(f)+t_\Phi(1-f)&\leq \lp \frac{1}{2}\rp^{d-1}+\sum_{\substack{k=2\\k \text{ even}}}^4 \frac{\alpha^k}{2^{d-k-1}} \sum_{I\in [d]^{(k)}} \E_{v\in G^r} \lp\prod_{i\in I} f_1(\phi_i(v))f_2(\phi_i(v))\rp + \alpha^6 \cdot 8^d.
\end{align*}
The key difference to \Cref{thm:ap4vector} is that we cannot split the directional and muting part. Instead we will prove the following bounds.
\begin{enumerate}
\item If $\vert I\vert =2$ and $p>4\cdot 5^4$, then there is a constant $C>0$, depending only on the matrix $M$ that induces $\Phi$, such that
\begin{align}\label{eq:cyclic-size-two-bound}
\lv \E_{w\in \Z_p^r} \prod_{i\in I} f_1(\phi_i(w))f_2(\phi_i(w)) \rv \leq C\frac{\log(p)+1}{\sqrt{p}}.
\end{align}
\item If $\vert I\vert=4$ and $I$ does not correspond to a 4-AP, then there is a constant $C>0$, depending only on $M$, such that if $p> C$, then
\begin{align}\label{eq:cyclic-non-4ap-bound}
\lv \E_{w\in \Z_p^r} \prod_{i\in I} f_1(\phi_i(w))f_2(\phi_i(w)) \rv \leq 80\beta^3 + C\frac{\log(p)+1}{\sqrt{p}}.
\end{align}
\item If $I$ corresponds to a 4-AP, then there is a constant $C>0$, again depending only on $M$, such that if $p> C$, then
\begin{align}\label{eq:cyclic-4ap-bound}
\E_{w\in \Z_p^r} \prod_{i\in I} f_2'(\phi_i(w)) \leq -2\beta^2\cdot \frac{2}{199^2}+C\frac{\log(p)+1}{\sqrt{p^n}}.
\end{align}
\end{enumerate}
Once these bounds are established, the proof concludes in exactly the same manner as that of \Cref{thm:ap4vector} under the condition that $p$ is sufficiently large.\\\\
We begin by showing \eqref{eq:cyclic-size-two-bound}. Let $I=\{i_1,i_2\}$. Reparametrizing the forms $\phi_{i_1},\phi_{i_2}$ such that $\phi_{i_j}'(w)=w_j$, we obtain
\begin{align*}
\lv \E_{w\in \Z_p^r} \prod_{i\in I}f_1(\phi_i(w))f_2(\phi_i(w))\rv &= \lv \sum_{\lambda\in \{\pm 1,\pm 3\}^2} \E_{x,y\in \Z_p} \beta^{\rho_2(\lambda)} \omega^{\lambda_1 (x^2+qx)+\lambda_2(y^2+qy)} f_1(x)f_1(y) \rv\\
&\leq \sum_{\lambda\in \{\pm 1,\pm 3\}^2} \lv \E_{x\in \Z_p} \omega^{\lambda_1 (x^2+qx)}f_1(x)\rv.
\end{align*}
But since $f_1$ is the weighted sum of two indicator functions of arithmetic progressions with common difference $1$ and $5$ and $p>4\cdot 5^4$, we can use \eqref{eq:ext-vanishing-phase-a} from \Cref{lem:phase-vanish-extended} to obtain \eqref{eq:cyclic-size-two-bound} with $C=96$, via the inequality
\begin{align*}
\sum_{\lambda\in \{\pm 1,\pm 3\}^2} \lv \E_{x\in \Z_p} \omega^{\lambda_1 (x^2+qx)}f_1(x)\rv \leq 16 \cdot 3 \cdot 2 \frac{\log(p)+1}{\sqrt{p}}=96\frac{\log(p)+1}{\sqrt{p}}.
\end{align*}
We now show \eqref{eq:cyclic-non-4ap-bound} for some set $I=\{i_1,\ldots,i_4\}$ that does not correspond to a 4-AP. First, we inspect the product of the muting parts. We have
\begin{align}\label{eq:cyclic-total-bound1}
\lv \E_{w\in \Z_p^r} \prod_{i\in I} f_1(\phi_i(w))f_2(\phi_i(w)) \rv&
\nonumber\\&\hspace*{-2cm} \leq\sum_{\lambda\in \{\pm 1,\pm 3\}^4}  \beta^{\rho_4(\lambda)}\lv \E_{w \in \Z_p^r} \lp \prod_{i\in I}f_1(\phi_i(w))\rp \omega^{\sum_{i\in I} \lambda_i(\phi_i(w)^2+q\phi_i(w))}\rv,
\end{align}
where once more $\rho_4(\lambda)=\lv\{ j \in [4]:\lambda_j=\pm 1\}\rv$. We will show that the exponent has a structure that will allow us to apply \Cref{lem:phase-vanish-extended}. But to be able to apply \Cref{lem:phase-vanish-extended}, we must keep track of the fraction bounds of the coefficients that occur. Let $C_1$ be an upper bound on the modulus of the entries of $M$, let $C_2=rC_1(C_1\sqrt{r})^r$ be the fraction bound from \Cref{lemma:reparametrization} for the coefficients of reparametrized systems, and let $C_3=12C_2^2$. \\\\
\textbf{Claim 1.} We may reparametrize $\{\phi_i\}_{i\in I}$ to $\{\phi_i'\}_{i\in I}$ such that the following holds. For all $\lambda\in \{\pm 1,\pm 3\}^4$ with $\rho_4(\lambda)\leq 2$ there is a permutation of coordinates
\begin{align*}
w_\lambda\colon \Z_p \times \Z_p \times \Z_p^{r-2} \rightarrow \Z_p^r \qquad (x,y,v)\mapsto w_\lambda(x,y,v),
\end{align*}
a map $g_\lambda\colon \Z_p^2\rightarrow \Z_p$ of the form
\begin{align*}
g_\lambda(x,y)=a_\lambda x^2+b_\lambda qx+d_\lambda xy,
\end{align*}
where $a_\lambda,b_\lambda,d_\lambda \in \Z$ are all $C_3$-fractions and at least one of $a_\lambda,b_\lambda,d_\lambda$ is coprime to $p$, and an additional map $h_\lambda\colon \Z_p\times \Z_p^{r-2}\rightarrow \Z_p$
such that
\begin{align*}
\forall x,y\in \Z_p, v\in \Z_p^{r-2}\colon \sum_{i\in I} \lambda_i(\phi_i'(w_\lambda(x,y,v))^2+q\phi_i'(w_\lambda(x,y,v)))=g_\lambda(x,y)+h_\lambda(y,v).
\end{align*}
\begin{proof}[Proof of Claim 1.]
As in the proof of \Cref{thm:ap4vector}, let $s$ be the dimension of the span of $\phi_{i_1},\ldots,\phi_{i_4}$. By pairwise linear independence of our forms, we know that $s\geq 2$. If $s=4$, all four forms are linearly independent and we reparametrize so that $\phi_{i_j}'(w)=w_j$. Letting $w_\lambda(x,y,v)=(x,y,v_1,\ldots,v_{r-2})$ then gives the claim with fraction bound $1<C_3$ on the coefficients.\\\\
If $s=3$, we may assume without loss of generality that the forms $\phi_{i_1},\phi_{i_2},\phi_{i_3}$ are independent. We reparametrize again so that $\phi_{i_j}'(w)=w_j$ for $j=1,2,3$. The fourth form $\phi'_{i_4}(w)$ is of the form $a_1w_1+a_2w_2+a_3w_3$, and by \Cref{lemma:reparametrization} the reparametrization can be carried out so that $a_1,a_2,a_3\in \Z$ are $C_2$-fractions. By the linear independence of pairs of forms, at least two of $a_1,a_2,a_3$ are coprime to $p$. If one of them, $a_3$ say, is zero mod $p$, then we let $w_\lambda(x,y,v)=(x,y,v_1,\ldots,v_{r-2})$ in order to obtain, for some $h_\lambda$, that for all $x,y\in \Z_p$ and $v\in \Z_p^{r-2}$
\begin{align*}
\sum_{i\in I} \lambda_i(\phi_i'(w_\lambda(x,y,v))^2+q\phi_i'(w_\lambda(x,y,v)))&=\\&\hspace*{-4.5cm}(\lambda_1+a_1^2\lambda_4)x^2+q(\lambda_1+a_1\lambda_4)x+2\lambda_4a_1a_2xy+h_\lambda(y,v).
\end{align*}
Note that $2\lambda_4a_1a_2$ is coprime to $p$ and that the coefficients of $x^2$, $x$ and $xy$ are indeed $C_3$-fractions.\\\\
If $a_1,a_2,a_3$ are all coprime to $p$, we reparametrize once more to 
$\phi_{i_j}''(w)=\phi_{i_j}'(Mw)$, where
\begin{align*}
M&=\begin{pmatrix}
A&0_{3\times (r-3)}\\0_{(r-3)\times 3}&I_{(r-3)\times (r-3)}
\end{pmatrix}&
A&=\begin{pmatrix}
1&0&0\\0&1&-a_2\inv\\0&0&a_3\inv
\end{pmatrix}
\end{align*}
and $a_2\inv,a_3\inv$ are representatives of the multiplicative inverses of $a_2,a_3$ mod $p$. It follows that $\phi_{i_1}''(w)=w_1$, $\phi_{i_2}''(w)=w_2-a_2\inv w_3$, $\phi_{i_3}''(w)=a_3\inv w_3$ and $\phi_{i_4}''(w)=a_1w_1+a_2w_2$. Because $w_1$ appears only in $\phi_{i_1}''$ and $\phi_{i_4}''$, we can proceed as before.\\\\
We are left with the case $s=2$. The same analysis as in the proof of \Cref{thm:ap4vector} shows that under the assumption that $\{\phi_{i}\}_{i\in I}$ does not form a 4-AP, the expression 
\begin{align*}
\sum_{i\in I} \lambda_i(\phi_i(w)^2+q\phi_i(w))
\end{align*}
cannot be zero for any $\lambda\in \{\pm 1, \pm 3\}$ with $\rho_4(\lambda)\leq 2$. We reparametrize in such a way that $\phi_{i_1}'(w)=w_1$ and $\phi_{i_2}'(w)=w_2$. This can be done while ensuring that $\phi_{i_3}'(w)=a_1w_1+b_1w_2$ and $\phi_{i_4}'(w)=a_2w_1+b_2w_2$ for $C_2$-fractions $a_1,a_2,b_1,b_2$. The expression
\begin{align*}
\sum_{i\in I} &\lambda_i(\phi_i'(w)^2+q\phi_i'(w))=(\lambda_1+a_1^2\lambda_3+a_2^2\lambda_4) w_1^2+(\lambda_2+b_1^2\lambda_3+b_2^2\lambda_4)w_2^2\\&\hspace*{1.5cm}+2(a_1b_1\lambda_3+a_2b_2\lambda_4)
w_1w_2
+(\lambda_1+a_1\lambda_3+a_2\lambda_4)qw_1 + (\lambda_2+b_1\lambda_3+b_2\lambda_4)qw_2.
\end{align*}
is non-zero. Therefore, one of the coefficients appearing in it must be non-zero. Setting $w_\lambda(x,y,v)$ to be $(x,y,v_1,\ldots,v_{r-2})$ or $(y,x,v_1,\ldots,v_{r-2})$ accordingly gives the claim.
\end{proof}
We reparametrize $\{\phi_i\}_{i\in I}$ according to the claim, but to keep the notation as simple as possible, we will continue to write $\phi_i$ instead of $\phi_i'$ for the reparametrized forms. While the coefficients of the $\phi_i$ may not be the original ones any longer, they are $C_3$-fractions.\\\\
There are $4\cdot 16+16=80$ different $\lambda\in \{\pm 1,\pm 3\}^4$ such that $\rho_4(\lambda)>2$ and since $f_1$ takes values in $[-1,1]$, we may bound \eqref{eq:cyclic-total-bound1} as
\begin{align}\label{eq:cyclic-total-bound2}
\lv \E_{w\in \Z_p^r} \prod_{i\in I} f_1(\phi_i(w))f_2(\phi_i(w)) \rv&\nonumber\\
&\hspace*{-2.5cm} \leq 80\beta^3 + \sum_{\substack{\lambda\in \{\pm 1,\pm 3\}^4\\\rho_4(\lambda)\leq 2}}  \E_{v \in \Z_p^{r-2}} \E_{y\in \Z_p} \lv  \E_{x\in \Z_p}\lp \prod_{i\in I}f_1(\phi_i(w_\lambda(x,y,v)))\rp \omega^{g_\lambda(x,y)}\rv.
\end{align}
Fix now some $\lambda$ with $\rho_4(\lambda)\leq 2$ and maps $w_\lambda, g_\lambda$ as in Claim 1. We also fix two specific elements $y\in \Z_p,v\in \Z_p^{r-2}$. Our directional part $f_1$ is of the form $f_1=1_{A_1}-2\cdot 1_{A_5}$ where $A_1$ and $A_5$ are projections of integer arithmetic progressions lying in $[p-1]_0$ with common difference 1 and 5, respectively. We now expand the product of muting parts, which yields
\begin{align}\label{eq:cyclic-direction-pre}
\lv \E_{x\in \Z_p} \omega^{g_\lambda(x,y)} \prod_{i\in I} f_1(\phi_i(w_\lambda(x,y,v))) \rv&\nonumber \\&\hspace*{-2.5cm}\leq 16 \sum_{\mu\in \{1,5\}^4} \lv \E_{x\in \Z_p} \omega^{g_\lambda(x,y)}\prod_{j=1}^4 1_{A_{\mu_j}}(\phi_{i_j}(w_\lambda(x,y,v))) \rv.
\end{align}
Observe that, because both the forms $\phi_i$ and the coordinate permutation $w_\lambda$ are linear, we have for all $x\in \Z_p$ that
\begin{align*}
1_{A_{\mu_j}}(\phi_{i_j}(w_\lambda(x,y,v)))=1_{A_{\mu_j}-\phi_{i_j}(w_\lambda(0,y,v))}(\phi_{i_j}(w_\lambda(x,0,0)))
\end{align*}
and $\phi_{i_j}(w_\lambda(x,0,0))=c_jx$ for some $C_3$-fraction $c_j$. If this coefficient is a multiple of $p$, then $1_{A_{\mu_j}}(\phi_{i_j}(w_\lambda(x,y,v)))$ does not depend on $x$ and we may bound \eqref{eq:cyclic-direction-pre} by pulling this factor out of the inner average and the modulus. If, on the other hand, $c_j$ is coprime to $p$, then there is a $c_j\inv\in \Z$ so that $c_j\inv c_j\equiv 1$ mod $p$. Let $B_{\mu,j}:=c_j\inv(A_{\mu_j}-\phi_{i_j}(w_\lambda(0,y,v)))$ and note that
\begin{align*}
1_{A_{\mu_j}-\phi_{i_j}(w_\lambda(0,y,v))}(\phi_{i_j}(w_\lambda(x,0,0)))=1_{B_{\mu,j}}(x).
\end{align*} 
Because $c_j$, and therefore also $c_j\inv$, are $C_3$-fractions, $B_{\mu,j}$ is a projection of an integer arithmetic progression whose common difference is a $5C_3$-fraction. Hence, by \Cref{lem:splitting} we can split $B_{\mu,j}$ into $m\leq 15C_3$ projections $B_{\mu,j}^1,\ldots,B_{\mu,j}^m$ of distinct integer APs,\footnote{Technically, the number of splitting progressions $m$ may vary with $\mu$ and $j$, but we simply take the maximal $m$ and add empty progressions where needed.} each with a common difference less then $5C_3$ and each lying in $[p-1]_0$. We let $J=\{j\in [4]:c_j\neq 0\}$ and define for each $\kappa\in [m]^J$
\begin{align*}
D_\mu^\kappa=\bigcap_{j\in J} B_{\mu,j}^{\kappa_j}
\end{align*}
to obtain
\begin{align}\label{eq:direction-bound1}
\lv \E_{x\in \Z_p} \omega^{g_\lambda(x,y)} \prod_{i\in I} f_1(\phi_i(w_\lambda(x,y,v))) \rv&\leq 16 \sum_{\mu\in \{1,5\}^4} \lv \E_{x\in \Z_p} \omega^{g_\lambda(x,y)} \prod_{j\in J} \lp \sum_{k=1}^m 1_{B_{\mu,j}^k}(x)\rp\rv \nonumber\\
&\leq 16 \sum_{\mu\in \{1,5\}^4} \sum_{\kappa\in [m]^J}\lv\E_{x\in \Z_p} \omega^{g_\lambda(x,y)}\prod_{j\in J} 1_{B_{\mu,j}^{\kappa_j}}(x)\rv\nonumber\\
&=16 \sum_{\mu\in \{1,5\}^4} \sum_{\kappa\in [m]^J} \lv\E_{x\in \Z_p} \omega^{g_\lambda(x,y)} \cdot 1_{D_\mu^\kappa}(x)\rv.
\end{align}
Because the $B_{\mu,j}^{\kappa_j}$ are projections of integer APs in $[p-1]_0$, the same is true of $D_\mu^\kappa$. What is more, the common difference of the pre-image in $\Z$ of $D_\mu^\kappa$ is at most the least common multiple of the common differences of the pre-images in $\Z$ of the progressions $B_{\mu,j}^{\kappa_j}$. This least common multiple can be bounded by $(5C_3)^{\lv J\rv}\leq (5C_3)^4=:C_4$. We can now apply \Cref{lem:phase-vanish-extended} with the fraction bound $C_4$ as long as $p>4C_4^4$. Indeed, recall that $g_\lambda(x,y)$ is of the form $a_\lambda x^2+b_\lambda qx+d_\lambda xy$, where $a_\lambda,b_\lambda,d_\lambda$ are $C_3$-fractions (and hence $C_4$-fractions) and not all are multiples of $p$. Keeping $v\in \Z_p^{r-2}$ fixed, we average \eqref{eq:direction-bound1} over $y$ to get
\begin{align}\label{eq:direction-bound2}
\E_{y\in \Z_p} \lv \E_{x\in \Z_p} \omega^{g_\lambda(x,y)} \prod_{i\in I} f_1(\phi_i(w_\lambda(x,y,v))) \rv&\leq 16 \sum_{\mu\in \{1,5\}^4} \sum_{\kappa\in [m]^J} \E_{y\in \Z_p}\lv\E_{x\in \Z_p} \omega^{g_\lambda(x,y)} \cdot 1_{D_\mu^\kappa}(x)\rv\nonumber\\
&\leq 16 \sum_{\mu\in \{1,5\}^4} \sum_{\kappa\in [m]^J} 8C_4^2 \frac{\log(p)+1}{\sqrt{p}}\nonumber\\
&\leq 256 \cdot (15C_3)^4 \cdot 8C_4^2 \frac{\log(p)+1}{\sqrt{p}}.
\end{align} 
Inserting this bound into \eqref{eq:cyclic-total-bound2} gives
\begin{align*}
\lv \E_{w\in \Z_p^r} \prod_{i\in I} f_1(\phi_i(w))f_2(\phi_i(w)) \rv&\leq 80\beta^3 + \sum_{\substack{\lambda\in \{\pm 1,\pm 3\}^4\\\rho_4(\lambda)\leq 2}}  256 \cdot (15C_3)^4 \cdot 8C_4^2 \frac{\log(p)+1}{\sqrt{p}}\\
&\leq 80\beta^3 + 176 \cdot 256 \cdot (15C_3)^4 \cdot 8C_4^2 \frac{\log(p)+1}{\sqrt{p}}.
\end{align*}
Letting $C:=\max(4C_4^4, 176\cdot 256 \cdot (15C_3)^4 \cdot 8C_4^2)$ yields \eqref{eq:cyclic-non-4ap-bound}.\\\\
Finally, we have to show \eqref{eq:cyclic-4ap-bound} assuming that $\{\phi_i\}_{i\in I}$ forms a 4-AP. By relabeling the indices in $I$ and reparametrizing, we can ensure that $\phi_{i_j}(w)=w_1+(j-1)w_2$. In particular, the forms do not depend on coordinates of $w$ other than $w_1$ and $w_2$. We write $x:=w_1$, $y:=w_2$ and $\phi_i(x,y)$ instead of $\phi_i(w)$. For $\lambda\in \{\pm 1,\pm 3\}$, let
\begin{align*}
g_\lambda\colon \Z_p^2\rightarrow \Z_p \qquad (x,y)\mapsto \sum_{i\in I} \lambda_i(\phi_i(x,y)^2+\phi_i(x,y)),
\end{align*}
so that
\begin{align}\label{eq:cyclic-ap-total-bound}
\E_{w \in \Z_p^r} \lp \prod_{i\in I}f_1(\phi_i(w))f_2(\phi_i(w))\rp&= \sum_{\lambda\in \{\pm 1,\pm 3\}^4} \beta^{\rho_4(\lambda)} \E_{x,y\in \Z_p} \omega^{g_\lambda(x,y)}\lp \prod_{i\in I}f_1(\phi_i(x,y))\rp.
\end{align}
If $g_\lambda$ is zero, which is certainly the case for $\lambda\in \{\pm (1,-3,3,-1)\}$, then by choice of $f_1$ we have that
\begin{align}\label{eq:cyclic-ap-muting-bound}
\E_{x,y\in \Z_p} \omega^{g_\lambda(x,y)}\lp \prod_{i\in I}f_1(\phi_i(x,y))\rp=\E_{x,y\in \Z_p} \lp \prod_{i\in I}f_1(\phi_i(x,y))\rp<-\frac{2}{199^2}.
\end{align}
For all $\lambda$ with $g_\lambda\neq 0$, the function $g_\lambda$ is an expression in $x^2,y^2,xy,x$ and $y$. Without loss of generality, we may assume that at least one coefficient of the terms involving $x$, i.e., $x^2,xy$ and $x$, is a non-zero integer. At the same time, all coefficients are bounded by the constant $C_3:=36$. Letting again $C_4:=(5C_3)^4$ and conducting the same analysis of the directional part as before, we obtain, analogously to \eqref{eq:direction-bound2}, the bound
\begin{align}\label{eq:cyclic-ap-direction-bound}
\lv \E_{x,y\in \Z_p} \omega^{g_\lambda(x,y)} \prod_{i\in I} f_1(\phi_i(x,y)) \rv& \leq 256 \cdot (15C_3)^4 \cdot 8C_4^2 \frac{\log(p)+1}{\sqrt{p}}.
\end{align} 
Applying the bounds \eqref{eq:cyclic-ap-muting-bound} and \eqref{eq:cyclic-ap-direction-bound} in \eqref{eq:cyclic-ap-total-bound} gives 
\begin{align*}
\E_{w \in \Z_p^r} \lp \prod_{i\in I}f_1(\phi_i(w))f_2(\phi_i(w))\rp&\leq -\frac{2}{199^2}\beta^2 + 254 \cdot 256 \cdot (15C_3)^4 \cdot 8C_4^2 \frac{\log(p)+1}{\sqrt{p}}.
\end{align*}
This yields \eqref{eq:cyclic-4ap-bound} for $C:=\max(C_4^4,254 \cdot 256 \cdot (15C_3)^4 \cdot 8C_4^2)$, completing the proof.
\end{proof}
We conclude this section by pointing out that the deduction of \Cref{cor:ap4cyclic} from \Cref{thm:ap4cyclic} is again a simple application of \Cref{lem:maps-to-sets}.
\section{Unmutable subconfigurations}\label{sec:cube-missing-corner}
At this point, one might hope to be able to apply the muting method to prove a compact characterization of common configurations. However, it turns out that there are some configurations that stubbornly resist any attempt to be muted. In any finite Abelian group $G$, consider, for instance, the additive quadruple, i.e., the configuration given by the system
\begin{align*}
\AQ\colon G^3 \rightarrow G^4 \qquad (x,h_1,h_2) \mapsto (x,x+h_1,x+h_2,x+h_1+h_2).
\end{align*}
For any $c\in \Z$ we have
\begin{align*}
cx-c(x+h_1)-c(x+h_2)+c(x+h_1+h_2)=0.
\end{align*}
Therefore, \textit{no} muting part consisting of phase functions with linear phases will actually mute the additive quadruple. This is all the more important because the contribution an additive quadruple makes to the arithmetic multiplicity is never negative, even if $f$ is allowed to take negative values. 
Indeed, for $f\colon G\rightarrow [-1,1]$, the arithmetic multiplicity of $t_\AQ(f)$ is
\begin{align}
\E_{x,h_1,h_2\in G} f(x)f(x+h_1)f(x+h_2)f(x+h_1+h_2)= \E_{h_2\in G} \lp \E_{x\in G} f(x)f(x+h_2)\rp^2 \geq 0.
\end{align}
The additive quadruple can be muted by a quadratic phase function, but this is only useful if our overall configuration contains a subconfiguration that is both sensitive to quadratic phases and can yield a negative arithmetic multiplicity.\\\\ 
Containing an additive quadruple as a subconfiguration does not only present an obstacle to the muting method but in fact has the potential to genuinely influence the arithmetic multiplicity of a configuration. The proof of \Cref{thm:cube-missing-corner} illustrates this point. 
\cubemissingcorner*
\begin{proof}
We fix a group $G$ and a map $g\colon G\rightarrow [0,1]$. Without loss of generality we may assume that $\alpha:=\E(g)\geq 1/2$, otherwise we could replace $g$ by $1-g$. Let further $f=g-\alpha$ and note that $f$ takes values in $[-\alpha,1-\alpha]\subseteq [-\alpha,\alpha]$. We rewrite the arithmetic multiplicities of $g$ and $1-g$ as
\begin{align}\label{eq:missing-cube-main}
t_\Phi(g)+t_\Phi(1-g)&=\E_{x,h_1,h_2,h_3\in G} \prod_{i=1}^7 \lp \alpha + f(\phi_i(x,h_1,h_2,h_3)\rp+\prod_{i=1}^7 \lp 1-\alpha - f(\phi_i(x,h_1,h_2,h_3))\rp\nonumber\\
&=\sum_{k=0}^7 \lp \alpha^{7-k}+(-1)^{k}(1-\alpha)^{7-k}\rp \sum_{I\in [7]^{(k)}} \E_{x,h_1,h_2,h_3} \prod_{i\in I} f(\phi_i(x,h_1,h_2,h_3)).
\end{align}
The summand corresponding to $k=7$ is always zero and the summand corresponding to $k=0$ is our main term of $\alpha^7+(1-\alpha)^7$, which is bounded from below by $(1/2)^6$. So it remains to be shown that the rest of the sum is non-negative. To see this, fix a $k\in \{1,\ldots,6\}$ as well as a set $I\in [7]^{(k)}$ and consider
\begin{align}\label{eq:cube-missing-sub-start}
\E_{x,h_1,h_2,h_3\in G} \prod_{i\in I} f(\phi_i(x,h_1,h_2,h_3)).
\end{align}
We will refer to this as the \textit{contribution} that the subconfiguration associated with $I$ makes to \eqref{eq:missing-cube-main}, ignoring for the moment the prefactor $\alpha^{7-k}+(-1)^{k}(1-\alpha)^{7-k}$. For $k=1,2,3$, it can easily be seen that \eqref{eq:cube-missing-sub-start} must be zero, because at least one of the variables $h_1,h_2,h_3$, will appear in exactly one of the forms $\phi_i$ for $i\in I$. This allows us to take an inner average of $f(\phi_i(x,h_1,h_2,h_3))$ running over just that variable. But this amounts to taking the average over $f$, which is zero.\\\\
For $k=4$, $(\phi_i)_{i\in I}$ is given by a vector in one of four classes. The first class $A_1$ contains all the subconfigurations $(\phi_i)_{i\in I}$ with $I\in [7]^{(4)}$ such that one of the forms $\phi_i$ is independent of the others in the sense that for all $x\in G$ the set $\{(\phi_j(v))_{j\in I\setminus \{i\}}: v\in G^4, \phi_i(v)=x\}$ is identical. Because $\E(f)=0$, the contribution of any subconfiguration in $A_1$ is zero. It is straightforward, if a bit cumbersome, to check that all remaining subconfigurations lie in one of the following classes.\footnote{We refer the reader once more to \cite{versteegen-4ap-uncommon-v2-arxiv}, an earlier version of this article, where the forms in $A_1$ are listed explicitly.}
\begin{align*}
A_2&=\{(x,x+h_1,x+h_2,x+h_1+h_2),(x,x+h_1,x+h_3,x+h_1+h_3),\\&\hspace*{1.0cm}(x,x+h_2,x+h_3,x+h_2+h_3)\},\\
A_3&=\{(x+h_1,x+h_2,x+h_1+h_3,x+h_2+h_3),(x+h_1,x+h_3,x+h_1+h_2,x+h_2+h_3),\\&\hspace*{1.0cm}
(x+h_2,x+h_3,x+h_1+h_2,x+h_1+h_3)\},\\
A_4&=\{(x,x+h_1+h_2,x+h_1+h_3,x+h_2+h_3)\}.
\end{align*}
Note that the contributions of the individual subconfigurations in each class are equal by symmetry. For $A_2$, we obtain the arithmetic multiplicity of additive quadruples, $t_{\AQ}(f)$. A simple calculation shows that the same is in fact true for $A_3$.\\\\
The class $A_4$ requires a bit of additional care. If the order of $G$ is odd, then the arithmetic multiplicity vanishes, but since we do not know anything about $G$ we instead bound the contribution using Cauchy-Schwarz as follows.
\begin{align*}
&\hspace*{-1.0cm}\lv \E_{x,h_1,h_2,h_3\in G} f(x)f(x+h_1+h_2)f(x+h_1+h_3)f(x+h_2+h_3) \rv\\
&=\lv \E_{h_2,h_3} \lp \E_x f(x)f(x+h_2+h_3)\rp \lp \E_y f(y+h_2)f(y+h_3)\rp\rv\\
&\leq \lp \E_{h_2,h_3} \lp \E_x f(x)f(x+h_2+h_3)\rp^2\rp^{1/2} \lp \E_{h_2,h_3}\lp \E_y f(y+h_2)f(y+h_3)\rp^2\rp^{1/2}\\
&=\lp \E_{h} \lp \E_x f(x)f(x+h)\rp^2\rp^{1/2} \lp \E_{h}\lp \E_y f(y)f(y+h)\rp^2\rp^{1/2},
\end{align*}
which is the same as $t_{\AQ}(f)$. In summary, the subconfigurations of size four taken together make a contribution of at least $5t_{\AQ}(f)$ to the arithmetic multiplicity $t_{\Phi}(g)$.\\\\
We now consider the subconfigurations for which $k=\lv I\rv=5$. Once more, we let $A_1$ be the class of subconfigurations that contain a form which is independent of the others. These subconfigurations again make no contribution to the arithmetic multiplicity. This just leaves the following three forms, which we subsume in the class $A_2$.
\begin{align*}
A_2&=\{(x,x+h_1,x+h_1+h_2,x+h_1+h_3,x+h_2+h_3),\\&\hspace*{1.0cm}
(x,x+h_2,x+h_1+h_2,x+h_1+h_3,x+h_2+h_3),\\&\hspace*{1.0cm}
(x,x+h_3,x+h_1+h_2,x+h_1+h_3,x+h_2+h_3)\}.
\end{align*}
The contribution from the elements of $A_2$ is not necessarily zero. However, we can bound it from above in absolute value by $(1-\alpha)t_{\AQ}(f)$. Indeed, note that
\begin{align*}
\{(x,x+h_1,&x+h_1+h_2,x+h_1+h_3,x+h_2+h_3)^T\in G^5:x,h_1,h_2,h_3\in G\}\\&=\{v\in G^5: v_1-2v_2+v_3+v_4-v_5=0\}.
\end{align*}
Therefore, using \eqref{eq:char-ortho-sum}, we see that
\begin{align}\label{eq:cube-missing-vtx-fourier-five}
&\hspace*{-1.0cm}\lv \E_{x,h_1,h_2,h_3\in G} f(x)f(x+h_1)f(x+h_1+h_2)f(x+h_1+h_3)f(x+h_2+h_3)\rv\nonumber\\&=\lv\sum_{\gamma\in \hat{G}} \E_{v\in G^5}f(v_1)f(v_2)f(v_3)f(v_4)f(v_5)\gamma(v_1-2v_2+v_3+v_4-v_5)\rv.
\end{align}
We can use the multiplicativity of characters to split the coordinates of the average, and then it follows from the definition of the Fourier transform that \eqref{eq:cube-missing-vtx-fourier-five} is equal to
\begin{align*}
\lv \sum_{\gamma\in \hat{G}} \hat{f}(\gamma)^3\hat{f}(\gamma\inv) \hat{f}(\gamma^{-2})\rv
\leq \sum_\gamma \vert \hat{f}(\gamma)\vert^4 \vert \hat{f}(\gamma^2)\vert
\leq \max_{\gamma\in \hat{G}} \vert \hat{f}(\gamma)\vert \sum_{\chi\in \hat{G}} \vert \hat{f}(\chi)\vert^4. 
\end{align*}
Another application of \eqref{eq:char-ortho-sum} shows that $\sum_\chi \lv \hat{f}(\chi)\rv^4=t_{\AQ}(f)$. Furthermore, we know that $\hat{f}(1)=0$ and for $\gamma\neq 1$ we have
\begin{align*}
\lv \hat{f}(\gamma)\rv = \lv \hat{(1-g)}(\gamma)\rv \leq \E_{x\in G} \lv 1- g(x)\rv\cdot \lv\gamma(x)\rv = \E_x (1-g(x))=1-\alpha.
\end{align*}
Finally, we need to consider the subconfigurations of size 6. There we have the classes
\begin{align*}
A_1&=\{(x,x+h_1,x+h_2,x+h_3,x+h_1+h_2,x+h_1+h_3),\\&\hspace*{1.0cm}
(x,x+h_1,x+h_2,x+h_3,x+h_1+h_2,x+h_2+h_3),\\&\hspace*{1.0cm}
(x,x+h_1,x+h_2,x+h_3,x+h_1+h_3,x+h_1+h_3)\},\\
A_2&=\{(x,x+h_1,x+h_2,x+h_1+h_2,x+h_1+h_3,x+h_2+h_3),\\&\hspace*{1.0cm}
(x,x+h_1,x+h_3,x+h_1+h_2,x+h_1+h_3,x+h_2+h_3),\\&\hspace*{1.0cm}
(x,x+h_2,x+h_3,x+h_1+h_2,x+h_1+h_3,x+h_2+h_3)\},\\
A_3&=\{(x+h_1,x+h_2,x+h_3,x+h_1+h_2,x+h_1+h_3,x+h_2+h_3)\}.
\end{align*}
The elements of $A_1$ make non-negative contributions. Indeed,
\begin{align}\label{eq:missing-cube-positive-sixer}
\E_{x,h_1,h_2,h_3\in G} &f(x)f(x+h_1)f(x+h_2)f(x+h_3)f(x+h_1+h_2)f(x+h_1+h_3)\nonumber\\
&=\E_{x,y,h_2,h_3} f(x)f(y)f(x+h_2)f(x+h_3)f(y+h_2)f(y+h_3)\nonumber\\
&=\E_{h_2,h_3}\lp \E_{x} f(x)f(x+h_2)f(x+h_3)\rp^2.
\end{align}
We call this configuration \textit{additive hextuple}, denote it by $\AH$ and refer to the quantity in \eqref{eq:missing-cube-positive-sixer} as $t_\AH(f)$ accordingly. The contribution from elements of $A_2$ can be bounded in modulus by $t_\AH(f)$ using the Cauchy-Schwarz inequality. Indeed, we reparametrize
\begin{align*}
\E_{x,h_1,h_2,h_3\in G} f(x)f(x+h_1)f(x+h_2)f(x+h_1+h_2)f(x+h_1+h_3)f(x+h_2+h_3)
\end{align*}
as
\begin{align*}
\E_{h_2,h_3} \lp \E_x f(x)f(x+h_2)f(x+h_2+h_3)\rp \lp \E_y f(y)f(y+h_2)f(y+h_3)\rp,
\end{align*}
which is bounded in modulus by
\begin{align*}
\lp \E_{h_2,h_3} \lp \E_x f(x)f(x+h_2)f(x+h_2+h_3)\rp^2 \rp^{1/2} \lp \E_{h_2,h_3} \lp \E_y f(y)f(y+h_2)f(y+h_3)\rp^2 \rp^{1/2},
\end{align*}
which equals $t_{\AH}(f)$. This means the subconfigurations from $A_1$ will always balance out any negative contributions that may arise from $A_2$. Lastly, we bound the contribution of the unique element of $A_3$ by $t_{\AQ}(f)$. Observe first that
\begin{align*}
&\hspace*{-1.0cm}\lv \E_{x,h_1,h_2,h_3\in G} f(x+h_1)f(x+h_2)f(x+h_3)f(x+h_1+h_2)f(x+h_1+h_3)f(x+h_2+h_3)\rv\\
&= \lv\E_{s,a,b,c} f(a)f(b)f(c)f(s-c)f(s-b)f(s-a)\rv.
\end{align*}
We use the symmetry in $a,b$ and $c$ of the above expression to bound it by
\begin{align*}
\E_{s} \lv \E_a f(a)f(s-a)\rv^3\leq \max_{s\in G} \lv \E_{a} f(a)f(s-a)\rv \E_{s} \lp\E_a f(a)f(s-a)\rp^2\leq \E(f^2)\cdot t_{\AQ}(f).
\end{align*}
To bound this further, note that by Parseval we have
\begin{align*}
\E(f^2)=\sum_{\gamma\in \hat{G}} \vert\hat{f}(\gamma)\vert^2 = \sum_{\gamma} \lv \hat{g}(\gamma)\rv^2-\lv\hat{g}(1)\rv^2=\E(g^2)-\E(g)^2\leq \alpha-\alpha^2.
\end{align*}
Taking all of these contributions together and recalling that $\alpha\geq 1/2$, \eqref{eq:missing-cube-main} becomes
\begin{align*}
t_\Phi(g)+t_\Phi(1-g)&=\sum_{k=0}^7 \lp \alpha^{7-k}+(-1)^{k}(1-\alpha)^{7-k}\rp \sum_{I\in [7]^{(k)}} \E_{x,h_1,h_2,h_3} \prod_{i\in I} f(\phi_i(x,h_1,h_2,h_3))\\
&\geq\frac{1}{64}+5(\alpha^3+(1-\alpha)^3)t_{\AQ}(f)-3(\alpha^2-(1-\alpha)^2)(1-\alpha)t_{\AQ}(f)\\&\hspace*{2.0cm}-\alpha(1-\alpha)t_{\AQ}(f)\\
&=\frac{1}{64}+t_{\AQ}(f)(22\alpha^2-25\alpha+8).
\end{align*}
But the polynomial $22\alpha^2-25\alpha+8$ has a global minimum of $79/88$ and since $t_{\AQ}(f)\geq 0$, we have $t_\Phi(g)+t_\Phi(1-g)\geq 1/64$ as desired.
\end{proof}
Unlike configurations given by a single equation in an odd number of variables, the cube with a missing vertex has no obvious property that indicates its commonness. This suggests that it might be difficult to even formulate a conjecture concerning a characterization of common configurations, even more so because it seems difficult to say how ``tight'' the above proof is. To keep the monochromatic arithmetic multiplicity above $1/64$ in the final argument, it would have sufficed if the configuration had only had three additive quadruples instead of six. It is conceivable that if one employed more sophisticated bounding methods on the negative contributions, even a single additive quadruple might have sufficed.\\\\ 
The following question may provide valuable insight into a possible characterization of common configurations. For the sake of simplicity, we will restrict its formulation to $\F_p^n$ for a fixed prime $p$. 
\begin{question}\label{qu:dominating-conf}
Let $M\in \Z^{d\times r}$ and denote for each $n\in \N$ by $\Phi_n$ the configuration induced by $M$ in $\F_p^n$. We say that $M$ is \emph{dominated by the additive quadruple} if there exists a constant $C>0$ such that for all $n\in \N$ and all functions $f\colon \F_p^n\rightarrow [-1,1]$ we have
\begin{align}\label{eq:AQ-dominated}
\lv t_{\Phi_n}(f)\rv \leq C t_{\AQ}(f).
\end{align}
Which configurations are dominated by the additive quadruple?
\end{question}
There are two classes of ``trivial'' examples that are \textit{not} dominated by the additive quadruple. Firstly, if $d<4$ we may simply take $f$ to be a sufficiently small constant function to see that \eqref{eq:AQ-dominated} is not satisfied. Secondly, if the squares of the forms in $\Phi$ are linearly dependent, we may construct $f$ as an appropriate quadratic phase function such that $t_{\Phi}(f)$ is large and $t_{\AQ}(f)$ converges to 0 as $n$ approaches infinity. In the other direction, we may employ the Cauchy-Schwarz inequality as in the proof of \Cref{thm:cube-missing-corner} to see that some configurations are dominated by the additive quadruple, but a more sophisticated approach may be necessary to prove this property for others.
\printbibliography
\end{document}